\documentclass[fleqn,12pt,a4paper]{article} 
\usepackage[left=1.5cm,right=1.5cm,top=2cm,bottom=2cm]{geometry} 

\usepackage{amsmath,amssymb,amsthm} 
\usepackage{mathrsfs} 
\usepackage{setspace} 
\usepackage{times} 
\usepackage{bm} 
\usepackage{enumitem} 
\usepackage{titlesec} 
\usepackage{fancyhdr} 
\usepackage{indentfirst} 

\usepackage{graphicx}  
\usepackage{caption}   
\usepackage{float}     

\titleformat{\section}{\large\bfseries}{\thesection.}{0.5em}{} 
\titleformat{\subsection}{\normalsize\bfseries}{\thesubsection.}{0.5em}{} 
\titleformat{\subsubsection}{\normalsize\itshape}{\thesubsubsection.}{0.5em}{} 

\begin{document}
	\title{Well-posedness of McKean-Vlasov generalized multivalued BSDEs}
	\author{
		Qi LIU\thanks{\rm Department of Mathematics, Shandong University, Shanda Road, Jinan City, Shandong Province. E-mail: 202311949@mail.sdu.edu.cn} \quad Yanbo CHEN\thanks{\rm Department of Mathematics, Shandong University, Shanda Road, Jinan City, Shandong Province. E-mail: 201911977@mail.sdu.edu.cn }
	}
	\date{March 11, 2026}  
	\maketitle
	
	\begin{abstract}
		\small This paper investigates McKean-Vlasov backward stochastic variational inequalities (BSVIs) whose generator depends on the joint law of the solution. We first establish the existence and uniqueness of the solution under globally Lipschitz and linear growth conditions. The analysis is then extended to the more general case of locally Lipschitz and non-linear growth conditions.
	\end{abstract}
	
	\section{Introduction}
	
	The foundational result on the existence and uniqueness of solutions to backward stochastic differential equations (BSDEs) was established by Pardoux and Peng~\cite{MR1037747}.  BSDEs are a fundamental tool in stochastic analysis, with widespread applications in mathematical finance, stochastic control, and partial differential equations. Pardoux and Peng~\cite{MR1037747} proved the uniqueness and existence of solutions to nonlinear BSDEs under uniformly Lipschitz conditions.
	
	Backward stochastic variational inequalities (BSVIs) extend BSDEs by incorporating  subdifferential operators, providing a powerful framework for modeling problems with state constraints and nonsmooth dynamics. This framework is particularly suitable for representing constrained stochastic optimization problems. Pardoux and R\u{a}\c{s}canu~\cite{MR1642656} introduced BSVIs where the multivalued term arises as the subdifferential of a convex function. By employing the Yosida approximation, they established the existence and uniqueness of solutions under linear growth and monotonicity conditions. The BSVI framework was further extended to infinite-dimensional Hilbert spaces in~\cite{MR1729473}. Maticiuc and R\u{a}\c{s}canu~\cite{MR2610326} investigated parabolic variational inequalities with multivalued Neumann--Dirichlet boundary conditions.
	
	Mean-field stochastic differential equations (MFSDEs), also known as McKean--Vlasov equations, can be traced back to the seminal works of Kac~\cite{MR84985} and McKean~\cite{MR221595}. These models aim to capture the collective dynamics of large systems of interacting particles by approximating microscopic interactions through their average effect. As such, they have been widely applied in statistical mechanics, quantum mechanics, and quantum chemistry. Building upon this foundation, Buckdahn et al.~\cite{MR2546754} introduced the framework of mean-field backward stochastic differential equations (MFBSDEs) using a fully probabilistic approach.
	
	A particular case of BSVI is the reflected BSDE (RBSDE). The constraint in RBSDE is enforced by a non-decreasing process, which corresponds to the case where the convex function in BSVI is an indicator function for the set defined by the obstacle. Under non-linear growth conditions, Boufoussi and Mouchtabih~\cite{MR4637498} constructed suitable sequences of approximating generators to establish the existence and uniqueness of solutions to RBSDEs by employing approximation and stability techniques.
	
	This paper studies a McKean-Vlasov BSVI by incorporating a dependency on the joint law of the solution into the generator of the BSVI presented in \cite{MR2610326}. Under the globally Lipschitz condition, our result reduces to the RBSDE studied in \cite{li2023general}. We then relax the assumptions imposed on the generator.
	
	In section 2, we address the case under the globally Lipschitz and linear growth conditions. The uniqueness is derived using Itô's formula and the Burkholder-Davis-Gundy (BDG) inequality, while the existence is established via Yosida approximation and contraction mapping.  In section 3, we study the locally Lipschitz and non-linear growth conditions. The existence is obtained by constructing a suitable sequence of generators that satisfy a global Lipschitz condition, then proving the convergence of the Cauchy sequences and passing to the limit.
	
	\section{BSVIs with globally Lipschitz and linear growth conditions}
	
	Let \(\{W_t : t \geq 0\}\) be a \(d\)-dimensional standard Brownian motion defined on some complete probability space \((\Omega, \mathcal{F}, \mathbb{P})\). We denote by \(\{\mathcal{F}_t : t \geq 0\}\) the natural filtration generated by \(\{W_t : t \geq 0\}\) and augmented by \(\mathcal{N}\), the set of \(\mathbb{P}\)-null events of \(\mathcal{F}\), $\mathcal{F}_t := \sigma \{W_r : 0 \leq r \leq t\} \vee \mathcal{N}.$  $T$ is a fixed time. Let \(\{A_t : t \geq 0\}\) be a continuous one-dimensional increasing progressively measurable stochastic process (p.m.s.p.) satisfying \(A_0 = 0\). $\xi$ is an $\mathbb{R}^k$-valued square-integrable  \(\mathcal{F}_T\)-measurable random variable. \(\mathcal{L}(Y_t, Z_t)\) represents the joint distribution of solution. $\delta_0$ is the Dirac measure with mass at 0 $\in \mathbb{R}^k \times \mathbb{R}^{k \times d }$.
	
	We investigate the existence and uniqueness of a solution \((Y, Z)\) to the following BSVI:
	\begin{equation}
		\begin{cases}
			\mathrm{d}Y_t + F(t, Y_t, Z_t,\mathcal{L}(Y_t, Z_t))\,\mathrm{d}t + G(t, Y_t)\,\mathrm{d}A_t \in \partial\varphi(Y_t)\,\mathrm{d}t + \partial\psi(Y_t)\,\mathrm{d}A_t + Z_t \,\mathrm{d}W_t, \\
			\qquad\qquad\qquad\qquad\qquad\qquad\qquad\qquad\qquad\qquad\qquad\qquad\qquad 0 \leq t \leq T, \\
			Y_T = \xi. 
		\end{cases}
	\end{equation}
	
	\subsection{Assumptions and results}
	Let \(\lambda, \mu \geq 0\).
	
	Let $\mathcal{H}_k^{\lambda, \mu} \subset L^2(\mathbb{R}_+ \times \Omega, e^{\lambda s + \mu A_s} \mathbf{1}_{[0, T]}(s) \, \mathrm{d}s \otimes \mathrm{d}\mathbb{P}; \mathbb{R}^k) $ be the Hilbert space of p.m.s.p. \(f : \Omega \times [0, \infty) \to \mathbb{R}^k\) such that
	\[
	\|f\|_{\mathcal{H}_k^{\lambda, \mu}} = \left[ \mathbb{E} \left( \int_0^T e^{\lambda s + \mu A_s} |f(s)|^2 \, \mathrm{d}s \right) \right]^{1/2} < \infty.
	\]
	Let $\widetilde{\mathcal{H}}_k^{\lambda, \mu} \subset L^2(\mathbb{R}_+ \times \Omega, e^{\lambda s + \mu A_s} \mathbf{1}_{[0, T]}(s) \, \mathrm{d}A_s \otimes \mathrm{d}\mathbb{P}; \mathbb{R}^k)$ be the Hilbert space of p.m.s.p. \(f : \Omega \times [0, \infty) \to \mathbb{R}^k\) such that
	\[
	\|f\|_{\widetilde{\mathcal{H}}_k^{\lambda, \mu}} = \left[ \mathbb{E} \left( \int_0^T e^{\lambda s + \mu A_s} |f(s)|^2 \, \mathrm{d}A_s \right) \right]^{1/2} < \infty.
	\]
	The notation \(\mathcal{S}_k^{\lambda, \mu}\) stands for the Banach space of p.m.s.p. \(f : \Omega \times [0, \infty) \to \mathbb{R}^k\) such that
	\[
	\|f\|_{\mathcal{S}_k^{\lambda, \mu}} = \left[ \mathbb{E} \left( \sup_{0 \leq t \leq T} e^{\lambda t + \mu A_t} |f(t)|^2 \right) \right]^{1/2} < \infty.
	\]
	$\mathcal{P}_2(\mathbb{R})$ is the space of Borel probability measures on $(\mathbb{R}, \mathcal{B}(\mathbb{R}))$ with finite second order moment.
	\[
	W_2(\mu, \nu) := \inf \left\{ \left( \int_{\mathbb{R}^{2}} |x - y|^2 \, \pi(dx, dy) \right)^{\frac{1}{2}} : \pi \in \mathcal{P}_2(\mathbb{R}^2), \ \pi(\cdot \times \mathbb{R}) = \mu, \ \pi(\mathbb{R} \times \cdot) = \nu \right\},
	\]where $\mu, \nu \in \mathcal{P}_2(\mathbb{R})$ and $\mathcal{B}(\mathbb{R})$ is the Borel $\sigma$-field over $\mathbb{R}$. 
	
	From the properties of the $W_2$-distance presented in  \cite{MR3752669}, we obtain:
	\begin{equation*}
		\begin{aligned}
			&W_2^2(\mathcal{L}(Y,Z), \mathcal{L}(Y',Z')) \leq (W_2(\mathcal{L}(Y,Z), \mathcal{L}(Y,Z'))+W_2(\mathcal{L}(Y,Z'), \mathcal{L}(Y',Z')))^2\\
			&\leq ((\mathbb{E}[\|Z - Z'\|^2])^\frac{1}{2}+(\mathbb{E}[|Y - Y'|^2])^\frac{1}{2})^2 \leq 2(\mathbb{E}[\|Z - Z'\|^2]+\mathbb{E}[|Y - Y'|^2]).
		\end{aligned}
	\end{equation*}
	
	With respect to BSVI (1), we formulate the following assumptions:
	
	\begin{itemize}[label=$\bullet$] 
		\item Assumption I\\ \(F : \Omega \times [0, \infty) \times \mathbb{R}^k \times \mathbb{R}^{k \times d }  \times \mathcal{P}_2(\mathbb{R}^k \times \mathbb{R}^{k \times d }) \to \mathbb{R}^k\). There exist \(L \geq 0\) and \(\eta : [0, \infty) \times \Omega \to [0, \infty)\) p.m.s.p. such that for all \(t \geq 0\), \(y, y' \in \mathbb{R}^k\), \(z, z' \in \mathbb{R}^{k \times d}\),\(\Pi,\Pi' \in \mathcal{P}_2(\mathbb{R}^k \times \mathbb{R}^{k \times d }) \),
		\begin{enumerate}[label=(\roman*)]
			\item \(F(\cdot, y, z,\Pi)\) is a p.m.s.p.
			\item \(y \mapsto F(\omega, t, y, z,\Pi) : \mathbb{R}^k \to \mathbb{R}^k\) is continuous a.s.
			\item \(|F(t, y,z,\Pi- F(t,y',z',\Pi')| \leq L(|y-y'|+ |z - z'|+W_2(\Pi,\Pi'))\) a.s.
			\item \(|F(t, y, z,\Pi)| \leq L(|y|+|z|+W_2(\Pi,\delta_0))+\eta_t \) a.s. 
		\end{enumerate}	
	\end{itemize}
	
	\begin{itemize}[label=$\bullet$] 
		\item  Assumption II\\ \(G : \Omega \times [0, \infty) \times \mathbb{R}^k \to \mathbb{R}^k\). There exist \(k_2\in \mathbb{R}, K \geq 0\) and \(\nu : [0, \infty) \times \Omega \to [0, \infty)\) is a p.m.s.p. such that for all \(t \geq 0\), \(y, y' \in \mathbb{R}^k\), \(z, z' \in \mathbb{R}^{k \times d}\),
		\begin{enumerate}[label=(\roman*)]
			\item \(G(\cdot, \cdot, y)\) is a p.m.s.p.
			\item \(y \mapsto G(\omega, t, y) : \mathbb{R}^k \to \mathbb{R}^k\) is continuous, a.s.
			\item \(|G(t, y) - G(t, y')|\leq K|y-y'|\) a.s.
			\item \(|G(t, y)| \leq K|y|+ \nu_t \), a.s. 
		\end{enumerate}
	\end{itemize}
	
	\begin{itemize}[label=$\bullet$] 
		\item  Assumption III\\ There exist $\lambda > 4+16 L^2$ and $\mu > 2+4 K^2$,
		\[
		\mathcal{M}(T)  :=\mathbb{E} \left[ e^{\lambda T + \mu A_T} \big(|\xi|^2 + \varphi(\xi) + \psi(\xi)\big) \right] 
		+ \mathbb{E} \int_0^T e^{\lambda s + \mu A_s} \big( |\eta_s|^2 \, \mathrm{d}s + |\nu_s|^2 \, \mathrm{d}A_s \big) < \infty.
		\]
	\end{itemize}
	
	\begin{itemize}[label=$\bullet$]
		\item Assumption IV\\We assume \(\varphi\) and \(\psi\) satisfy the following:
		\begin{enumerate}[label=(\roman*)]
			\item \(\varphi, \psi : \mathbb{R}^k \to (-\infty, +\infty]\) are proper convex lower semicontinuous (l.s.c) functions,
			\item \(\varphi(y) \geq \varphi(0) = 0,\quad \psi(y) \geq \psi(0) = 0.\) 
		\end{enumerate}
		
		The subdifferential of $\varphi(x)$ is defined by
		\[
		\partial \varphi(x) := \left\{ v \in \mathbb{R}^k : \langle v, y - x \rangle + \varphi(x) \leq \varphi(y), \, \forall y \in \mathbb{R}^k \right\}
		\]
		and similarly for \(\psi\).
		
		For \(\varepsilon > 0\), we define the convex \(C^1\)-function \(\varphi_\varepsilon\) via
		\[
		\varphi_\varepsilon(y) := \inf \left\{ \frac{1}{2\varepsilon} |v - y|^2 + \varphi(v) : v \in \mathbb{R}^k \right\},
		\]
		and similarly for \(\psi_\varepsilon\).
		
		From section 7.1 in \cite{MR2759829}, we  define $J_\varepsilon (y) := (I + \varepsilon \partial \varphi)^{-1} y$, which is the only minimizer of $\varphi_\varepsilon$. We substitute it into the definition and calculate the gradient $$\varphi_\varepsilon(y) = \varphi(J_\varepsilon(y)) + \dfrac{1}{2\varepsilon} |y - J_\varepsilon(y)|^2,\quad \nabla \varphi_\varepsilon(y) = \frac{y - J_\varepsilon (y)}{\varepsilon}.$$ 
		Similarly, $$\hat{J}_\varepsilon (y) := (I + \varepsilon \partial \psi)^{-1} y,\quad \psi_\varepsilon(y) = \frac{1}{2\varepsilon} |y - \hat{J}_\varepsilon (y)|^2 + \psi(\hat{J}_\varepsilon (y)),\quad \nabla \psi_\varepsilon(y) = \dfrac{y - \hat{J}_\varepsilon (y)}{\varepsilon} .$$
		
		As shown in \cite{MR274683}, \(\nabla \varphi_\varepsilon(y)\) and \(\nabla \psi_\varepsilon(y)\) are monotone Lipschitz functions.
		
		\item  Assumption V\\ Compatibility assumptions are introduced: $\text{for all } \varepsilon > 0, y \in \mathbb{R}^k \text{ and } z \in \mathbb{R}^{k \times d},$
		\begin{enumerate}[label=(\roman*)]
			\item \(\langle \nabla \varphi_\varepsilon(y), \nabla \psi_\varepsilon(y) \rangle \geq 0,\)
			\item \(\langle \nabla \varphi_\varepsilon(y), G(t, y) \rangle \leq \langle \nabla \psi_\varepsilon(y), G(t, y) \rangle^+,\)
			\item \(\langle \nabla \psi_\varepsilon(y), F(t, y, z,\mu) \rangle \leq \langle \nabla \varphi_\varepsilon(y), F(t, y, z,\mu) \rangle^+.\) 
		\end{enumerate}
	\end{itemize}
	
	\noindent\textbf{Definition 1}:  \((Y, Z, U, V)\) is a solution of BSVI (1) if:
	\begin{enumerate}
		\item[(a)] \(Y \in \mathcal{S}_k^{\lambda, \mu} \cap \mathcal{H}_k^{\lambda, \mu} \cap \widetilde{\mathcal{H}}_k^{\lambda, \mu}\), \quad \(Z \in \mathcal{H}_{k \times d}^{\lambda, \mu}\).
		\item[(b)] \(U \in \mathcal{H}_k^{\lambda, \mu}\),  \quad \(V \in \widetilde{\mathcal{H}}_k^{\lambda, \mu}\).
		\item[(c)] \(\mathbb{E} \int_0^T e^{\lambda s + \mu A_s} \left( \varphi(Y_s) \, ds + \psi(Y_s) \, dA_s \right) < \infty.\)
		\item[(d)] \((Y_t, U_t) \in \partial \varphi,\) \(d\mathbb{P} \times dt\), \quad   \((Y_t, V_t) \in \partial \psi,\)  \(d\mathbb{P} \times dA_t\)  \quad a.s. on \([0, T]\), i.e
        $$\langle U_t, u_t - Y_t \rangle + \varphi(Y_t) \leq \varphi(u_t), \, \quad \forall u \in  \mathcal{H}_k^{\lambda, \mu},\quad \langle V_t, v_t - Y_t \rangle + \varphi(Y_t) \leq \varphi(v_t), \, \quad \forall v \in  \tilde{ \mathcal{H}}_k^{\lambda, \mu},$$
		\item[(e)] The following equation is satisfied for all \(t \in [0,T]\)
		\begin{equation}
			\begin{aligned}
				&Y_t + \int_t^T U_s \, ds + \int_t^T V_s \, dA_s  = \xi + \int_t^T F(s, Y_s, Z_s,\mathcal{L}(Y_t, Z_t))ds \\&\quad+ \int_t^T G(s, Y_s) \, dA_s - \int_t^T Z_s \, dW_s.
			\end{aligned}
		\end{equation} 
	\end{enumerate}
	
	In the sequel, constant \( C \) depends on the parameters \( \lambda, \mu, K, L, T \), and its value may vary from line to line.
	
	\noindent\textbf{Proposition 2}: Let Assumptions I-IV be satisfied. If \((Y, Z, U, V)\) and \((\tilde{Y}, \tilde{Z}, \tilde{U}, \tilde{V})\) are solutions satisfying Definition 1 with terminal conditions \(\xi\) and \(\tilde{\xi}\), respectively. Then
	\begin{equation}
		\begin{aligned}
			&\mathbb{E} \int_0^T e^{\lambda s + \mu A_s} \left[ |Y_s - \tilde{Y}_s|^2 (\mathrm{d}s + \mathrm{d}A_s) + \|Z_s - \tilde{Z}_s\|^2 \mathrm{d}s \right]
			+ \mathbb{E} \sup_{0 \leq s \leq T} e^{\lambda s + \mu A_s} |Y_s - \tilde{Y}_s|^2\\& \leq C \, \mathbb{E} \,[ e^{\lambda T + \mu A_T} |\xi - \tilde{\xi}|^2]. 
		\end{aligned}
	\end{equation}
	
	\textbf{Proof.} From Itô's formula, we have
	\begin{equation*}
		\begin{aligned}
			&e^{\lambda t + \mu A_t} |Y_t - \tilde{Y}_t|^2 
			+ \int_t^T e^{\lambda s + \mu A_s} |Y_s - \tilde{Y}_s|^2 (\lambda \mathrm{d}s + \mu \mathrm{d}A_s)+ \int_t^T e^{\lambda s + \mu A_s} \|Z_s - \tilde{Z}_s\|^2 \mathrm{d}s \\
			&\quad+ 2 \int_t^T e^{\lambda s + \mu A_s} \langle Y_s - \tilde{Y}_s, U_s - \tilde{U}_s\rangle  \mathrm{d}s 
			+ 2\int_t^T e^{\lambda s + \mu A_s}  \langle Y_s - \tilde{Y}_s, V_s - \tilde{V}_s\rangle  \mathrm{d}A_s \\
			&= e^{\lambda T+ \mu A_T} |\xi - \tilde{\xi}|^2 \\&\quad+ 2 \int_t^T e^{\lambda s + \mu A_s}  \langle Y_s - \tilde{Y}_s, F(s, Y_s, Z_s,\mathcal{L}(Y_s, Z_s)) - F(s, \tilde{Y}_s, \tilde{Z}_s,\mathcal{L}(\tilde{Y}_s, \tilde{Z}_s)) \rangle \mathrm{d}s \\
			&\quad+ 2 \int_t^T e^{\lambda s + \mu A_s}  \langle Y_s - \tilde{Y}_s, G(s, Y_s) - G(s, \tilde{Y}_s)\rangle  \mathrm{d}A_s\\&\quad - 2 \int_t^T e^{\lambda s + \mu A_s}  \langle Y_s - \tilde{Y}_s, (Z_s - \tilde{Z}_s) \mathrm{d}W_s\rangle.\\
		\end{aligned}
	\end{equation*}
	For any $p > 0$, $2ab \leq p a^2 + \frac{1}{p} b^2$. We take $p = 4L^2$ here, then
	$$2L|Y_s - \tilde{Y}_s|\|Z_s - \tilde{Z}_s\|\leq 4L^2 |Y_s - \tilde{Y}_s|^2 + \frac{1}{4} \|Z_s - \tilde{Z}_s\|^2.$$
	With $W_2(\mathcal{L}(Y_t,Z_t),\mathcal{L}(Y'_t,Z'_t))\leq\sqrt{\mathbb{E}[\|Z_s - \tilde{Z}_s\|^2]}+\sqrt{\mathbb{E}[|Y_s - \tilde{Y}_s|^2]}$ and Assumptions I, II, we get:
	\begin{equation}
		\begin{aligned}
			&2 \langle Y_s - \tilde{Y}_s,F(s, Y_s, Z_s,\mathcal{L}(Y_s, Z_s)) - F(s, \tilde{Y}_s, \tilde{Z}_s,\mathcal{L}(\tilde{Y}_s, \tilde{Z}_s))\rangle \\
			&\leq 2 \langle |Y_s - \tilde{Y}_s|,L(|Y_s - \tilde{Y}_s|+\|Z_s - \tilde{Z}_s\|+\sqrt{\mathbb{E}[\|Z_s - \tilde{Z}_s\|^2]}+\sqrt{\mathbb{E}[|Y_s - \tilde{Y}_s|^2]})\rangle \\
			& \leq (2+9L^2)|Y_s - \tilde{Y}_s|^2+ \frac{1}{4}\|Z_s - \tilde{Z}_s\|^2+ \dfrac{1}{4}\mathbb{E}[\|Z_s - \tilde{Z}_s\|^2]  + L^2\mathbb{E}[|Y_s - \tilde{Y}_s|^2],\qquad\qquad
		\end{aligned}
	\end{equation}
	\begin{equation}
		\begin{aligned}
			&2 \langle Y_s - \tilde{Y}_s, G(s, Y_s) - G(s, \tilde{Y}_s) \rangle \leq  |Y_s - \tilde{Y}_s|^2+| G(s, Y_s) - G(s, \tilde{Y}_s)|^2\\ &\leq  (K^2+1)|Y_s - \tilde{Y}_s|^2.
		\end{aligned}
	\end{equation}
	Considering the definition 1 and the definition of subdifferentials, we have:
	\begin{equation}
		\langle Y_s - \tilde{Y}_s, U_s - \tilde{U}_s \rangle \geq 0, \quad 
		\langle Y_s - \tilde{Y}_s, V_s - \tilde{V}_s \rangle \geq 0.
	\end{equation}
	Therefore,
	\begin{equation}
		\begin{aligned}
			&e^{\lambda t+ \mu A_t} |Y_t - \tilde{Y}_t|^2 
			+ \int_t^T e^{\lambda s + \mu A_s} |Y_s - \tilde{Y}_s|^2 (\lambda \mathrm{d}s + \mu \mathrm{d}A_s)+\int_t^T e^{\lambda s + \mu A_s} \|Z_s - \tilde{Z}_s\|^2 \mathrm{d}s\\
			&\leq e^{\lambda T + \mu A_T} |\xi - \tilde{\xi}|^2 
			+ \int_t^T e^{\lambda s + \mu A_s}   (K^2+1)|Y_s - \tilde{Y}_s|^2  \mathrm{d}A_s \\
			&\quad +  \int_t^T e^{\lambda s + \mu A_s} \big( (2+9L^2)|Y_s - \tilde{Y}_s|^2+ \frac{1}{4}\|Z_s - \tilde{Z}_s\|^2\big) \mathrm{d}s\\
			&\quad +  \int_t^T e^{\lambda s + \mu A_s} \big(\dfrac{1}{4}\mathbb{E}[\|Z_s - \tilde{Z}_s\|^2]  + L^2\mathbb{E}[|Y_s - \tilde{Y}_s|^2]\big) \mathrm{d}s\\
			&\quad - 2 \int_t^T e^{\lambda s + \mu A_s}  \langle Y_s - \tilde{Y}_s, (Z_s - \tilde{Z}_s) \mathrm{d}W_s\rangle.\\
		\end{aligned}
	\end{equation}
	We get the following by using BDG inequality  
	\begin{equation}
		\begin{aligned}
			&\mathbb{E}\sup_{0 \leq r \leq T}\big|\int_r^T e^{\lambda s + \mu A_s}  \langle Y_s - \tilde{Y}_s, (Z_s - \tilde{Z}_s)\mathrm{d}W_s\rangle\big|\\
			&\leq C \mathbb{E} \left( \int_0^T e^{2(\lambda s + \mu A_s)} |Y_s - \tilde{Y}_s|^2\|Z_s - \tilde{Z}_s\|^2 ds \right)^{1/2}\\
			& \leq \frac{1}{4} \mathbb{E} \sup_{0 \leq s \leq T} e^{\lambda s + \mu A_s} |Y_s - \tilde{Y}_s|^2 + C \mathbb{E} \int_0^T e^{\lambda s + \mu A_s} \|Z_s - \tilde{Z}_s\|^2 ds.\qquad\qquad\qquad\qquad
		\end{aligned}
	\end{equation}
	By using Assumption III, taking expectation on (7) and setting $t=0$, we obtain 
	\begin{equation*}
		\mathbb{E} \int_0^T e^{\lambda s + \mu A_s} \left[ |Y_s - \tilde{Y}_s|^2 (\mathrm{d}s + \mathrm{d}A_s) + \|Z_s - \tilde{Z}_s\|^2 \mathrm{d}s \right]\leq C \, \mathbb{E} \,[ e^{\lambda T + \mu A_T} |\xi - \tilde{\xi}|^2]. 
	\end{equation*}
	By taking the supremum over time and the expectation in (7) and combining the result with (8), we obtain
	\begin{equation*}
		\mathbb{E} \sup_{0 \leq s \leq T} e^{\lambda s + \mu A_s} |Y_s - \tilde{Y}_s|^2 \leq C \, \mathbb{E} \,[ e^{\lambda T + \mu A_T} |\xi - \tilde{\xi}|^2].
	\end{equation*}
	The desired estimate (3) is a direct consequence of combining the preceding results.

	We now present the main result of this section.
	
	\noindent\textbf{Theorem 3}: Let Assumptions I-V be satisfied. There exists a unique solution \((Y, Z, U, V)\) for BSVI (1).
	
	\subsection{Existence and uniqueness}
	
	Consider the approximating equation
	\begin{equation}
		\begin{aligned}
			&Y_t^\varepsilon + \int_t^T \nabla \varphi_\varepsilon(Y_s^\varepsilon) \mathrm{d}s + \int_t^T \nabla \psi_\varepsilon(Y_s^\varepsilon) \mathrm{d}A_s 
			= \xi + \int_t^T F(s, Y_s^\varepsilon, Z_s^\varepsilon,\mathcal{L}(Y_s^\varepsilon, Z_s^\varepsilon)) \mathrm{d}s \\
			&\quad+ \int_t^T G(s, Y_s^\varepsilon) \mathrm{d}A_s 
			- \int_t^T Z_s^\varepsilon \mathrm{d}W_s, \quad t \in [0,T].
		\end{aligned}
	\end{equation}
	
	\noindent\textbf{Proposition 4}: Let Assumptions I-V be satisfied. Then there exists a unique solution  $(Y, Z, U, V) \in (\mathcal{S}_k^{\lambda, \mu} \cap \mathcal{H}_k^{\lambda, \mu} \cap \widetilde{\mathcal{H}}_k^{\lambda, \mu}) \times \mathcal{H}_{k \times d}^{\lambda, \mu} \times \mathcal{H}_k^{\lambda, \mu}\times \tilde{\mathcal{H}}_k^{\lambda, \mu}$ for (9).
	
	\textbf{Proof.} We define $\varPhi(Y,Z)=(\hat{Y},\hat{Z})$ as following and prove $\varPhi$ is a contraction mapping.  
	\begin{equation*}
		\begin{aligned}
			&\hat{Y}_t^\varepsilon + \int_t^T \nabla \varphi_\varepsilon(\hat{Y}_s^\varepsilon) \mathrm{d}s + \int_t^T \nabla \psi_\varepsilon(\hat{Y}_s^\varepsilon) \mathrm{d}A_s 
			= \xi + \int_t^T F(s, Y_s^\varepsilon, Z_s^\varepsilon,\mathcal{L}(Y_s^\varepsilon, Z_s^\varepsilon)) \mathrm{d}s \\
			&\quad+ \int_t^T G(s, Y_s^\varepsilon) \mathrm{d}A_s 
			- \int_t^T \hat{Z}_s^\varepsilon \mathrm{d}W_s.
		\end{aligned}
	\end{equation*}
	Since $\varepsilon$ is fixed, we omit it from the notation for simplicity during the proof.
	
	Similar to (4) and (5), we have
	\begin{equation*}
		\begin{aligned}
			& 2 \langle  \hat{Y}_s^{1} - \hat{Y}_s^{2}, F(s, Y_s^{1}, Z_s^{1},\mathcal{L}(Y_s^{1}, Z_s^{1})) - F(s,Y_s^{2}, Z_s^{2}\mathcal{L}(Y_s^{2}, Z_s^{2})) \rangle\\
			&\leq  16 L^2|\hat{Y}_s^{1} - \hat{Y}_s^{2}|^2 + \frac{1}{4 }|Y_s^{1} - Y_s^{2}|^2 + \frac{1}{4 }\|Z_s^{1} - Z_s^{2}\|^2 + \dfrac{1}{4 }\mathbb{E}\|Z_s^{1} - Z_s^{2}\|^2  + \frac{1}{4 }\mathbb{E}|Y_s^{1} - Y_s^{2}|^2,\\
			&2\langle  \hat{Y}_s^{1} - \hat{Y}_s^{2}, G(s, Y_s^{1}) - G(s, Y_s^{2})\rangle \leq 4 K^2|\hat{Y}_s^{1} - \hat{Y}_s^{2}|^2+\frac{1}{4 }|Y_s^{1}-Y_s^{2}|^2.
		\end{aligned}
	\end{equation*}
	Using Itô's formula and the above inequalities, we get
	\begin{equation*}
		\begin{aligned}
			&e^{\lambda t + \mu A_t} |\hat{Y}_t^{1} - \hat{Y}_t^{2} |^2 
			+ \int_t^T e^{\lambda s + \mu A_s} |\hat{Y}_s^{1} - \hat{Y}_s^{2} |^2 (\lambda \mathrm{d}s + \mu \mathrm{d}A_s) +\int_t^T e^{\lambda s + \mu A_s} \|\hat{Z}_s^{1} - \hat{Z}_s^{2} \|^2 \mathrm{d}s\\
			&\leq  \int_t^T e^{\lambda s + \mu A_s} \big(16 L^2 |\hat{Y}_s^{1} - \hat{Y}_s^{2}|^2+ \frac{1}{4}|Y_s^{1} - Y_s^{2}|^2 + \frac{1}{4 }\|Z_s^{1} - Z_s^{2}\|^2\big) \mathrm{d}s \\
			&\quad+  \int_t^T e^{\lambda s + \mu A_s} \big(\dfrac{1}{4 }\mathbb{E}\|Z_s^{1} - Z_s^{2}\|^2 + \frac{1}{4}\mathbb{E}|Y_s^{1} - Y_s^{2}|^2\big)\mathrm{d}s \\
			&\quad+ \int_t^T e^{\lambda s + \mu A_s} \big(4K^2|\hat{Y}_s^{1} - \hat{Y}_s^{2}|^2+\frac{1}{4 }|Y_s^{1}-Y_s^{2}|^2\big)\mathrm{d}A_s.\\
			&\quad - 2 \int_t^T e^{\lambda s + \mu A_s}  \langle  \hat{Y}_s^{1} - \hat{Y}_s^{2}, ( \hat{Z}_s^{1} - \hat{Z}_s^{2}) \mathrm{d}W_s\rangle. \\
		\end{aligned}
	\end{equation*}
	By setting $t=0$ and taking expectations, and invoking Assumption III, we obtain
	\begin{equation*}
		\begin{aligned}
			&\mathbb{E}\int_0^Te^{\lambda s + \mu A_s} |\hat{Y}_s^{1} - \hat{Y}_s^{2} |^2 (\mathrm{d}s +\mathrm{d}A_s)+ \mathbb{E} \int_0^T e^{\lambda s + \mu A_s} \|\hat{Z}_s^{1} - \hat{Z}_s^{2} \|^2 \mathrm{d}s \\
			&\leq \frac{1}{2 }\big(\mathbb{E} \int_0^T e^{\lambda s + \mu A_s} |Y_s^{1} - Y_s^{2}|^2(\mathrm{d}s+\mathrm{d}A_s)+\mathbb{E} \int_0^T e^{\lambda s + \mu A_s} \|Z_s^{1} - Z_s^{2} \|^2 \mathrm{d}s\big).
		\end{aligned}
	\end{equation*}
	Therefore, $\varPhi$ is a contraction mapping. Equation (9) admits a unique solution in $(\mathcal{H}_k^{\lambda, \mu} \cap \widetilde{\mathcal{H}}_k^{\lambda, \mu}) \times \mathcal{H}_{k \times d}^{\lambda, \mu}$. 
	
	By using Itô's formula for \(e^{\lambda t + \mu A_t} |Y_t^\varepsilon|^2\), we have
	\begin{equation*}
		\begin{aligned}
			&e^{\lambda t + \mu A_t} |Y_t^\varepsilon|^2 
			+\int_t^T e^{\lambda s + \mu A_s} \big(|Y_s^\varepsilon|^2 (\lambda \mathrm{d}s + \mu \mathrm{d}A_s) 
			+ \|Z_s^\varepsilon\|^2 \mathrm{d}s\big) \\
			&\quad + 2 \int_t^T e^{\lambda s + \mu A_s} \big(\langle Y_s^\varepsilon, \nabla \varphi_\varepsilon(Y_s^\varepsilon) \rangle \mathrm{d}s
			+ \langle Y_s^\varepsilon, \nabla \varphi_\varepsilon(Y_s^\varepsilon) \rangle \mu \mathrm{d}A_s \big) \\
			&= e^{\lambda T + \mu A_T} |\xi|^2 - 2 \int_t^T e^{\lambda s + \mu A_s} \langle Y_s^\varepsilon, Z_s^\varepsilon \mathrm{d}W_s \rangle\\
			&\quad + 2 \int_t^T e^{\lambda s + \mu A_s} \big(\langle Y_s^\varepsilon, F(s, Y_s^\varepsilon, Z_s^\varepsilon,\mathcal{L}(Y_s^\varepsilon, Z_s^\varepsilon)) \rangle \mathrm{d}s
			+ \langle Y_s^\varepsilon, G(s, Y_s^\varepsilon) \rangle \mathrm{d}A_s \big).
		\end{aligned}
	\end{equation*}
	From Assumptions I and II, we have
	\begin{equation*}
		\begin{aligned}
			&2 \langle Y_s^\varepsilon, F(s, Y_s^\varepsilon, Z_s^\varepsilon,\mathcal{L}(Y_s^\varepsilon, Z_s^\varepsilon)) \rangle\\&\leq (3+9L^2)|Y_s^\varepsilon|^2+ \frac{1}{4}\|Z_s^\varepsilon\|^2+ \dfrac{1}{4}\mathbb{E}[\|Z_s^\varepsilon\|^2]  + L^2\mathbb{E}[|Y_s^\varepsilon|^2]  + \eta_s^2,\\
			&2 \langle Y_s^\varepsilon, G(s, Y_s^\varepsilon) \rangle =2 \langle Y_s^\varepsilon, G(s, Y_s^\varepsilon)-G(s,0) \rangle + 2 \langle Y_s^\varepsilon, G(s,0) \rangle \\&\leq(K^2+2)|Y_s^\varepsilon|^2 + \nu_s^2. \qquad\qquad\qquad\qquad
		\end{aligned}
	\end{equation*}
	Taking the supremum over time and the expectation, and applying (8), we obtain the following result under Assumption IV:
	\begin{equation}
		\mathbb{E} \sup_{0 \leq s \leq T} e^{\lambda s + \mu A_s} | Y_s^\epsilon |^2 \leq C \mathcal{M}(T).
	\end{equation}
	Therefore, equation (9) admits a unique solution in $ (\mathcal{S}_k^{\lambda, \mu} \cap \mathcal{H}_k^{\lambda, \mu} \cap \widetilde{\mathcal{H}}_k^{\lambda, \mu}) \times \mathcal{H}_{k \times d}^{\lambda, \mu}$.

	\subsubsection{A priori estimate}
	
	\noindent\textbf{Proposition 5}: Let Assumptions I-IV be satisfied. Then
	\begin{equation}
		\begin{aligned}
			&\mathbb{E} \left[ \sup_{0 \leq s \leq T} e^{\lambda s + \mu A_s} |Y_s^\varepsilon|^2 
			+ \int_0^T e^{\lambda s + \mu A_s} \big(|Y_s^\varepsilon|^2 + \|Z_s^\varepsilon\|^2\big) \mathrm{d}s
			+ \int_0^T e^{\lambda s + \mu A_s} |Y_s^\varepsilon|^2 \mathrm{d}A_s \right]\\&\leq C \mathcal{M}(T).
		\end{aligned}
	\end{equation}
	
	\textbf{Proof.} Following a similar argument to the proof of Proposition 4, we obtain
	\begin{equation}
		\begin{aligned}
			&e^{\lambda t + \mu A_t} |Y_t^\varepsilon|^2 + \int_t^T e^{\lambda s + \mu A_s} |Y_s^\varepsilon|^2 ((\lambda-3-9L^2)  \mathrm{d}s + (\mu-2-K^2)\mathrm{d}A_s) \\
			&\quad + \int_t^T e^{\lambda s + \mu A_s} \dfrac{3}{4}\|Z_s^\varepsilon\|^2 \mathrm{d}s\\
			&\leq e^{\lambda T+ \mu A_T} |\xi|^2 + \int_t^T e^{\lambda s + \mu A_s} \big( \dfrac{1}{4}\mathbb{E}[\|Z_s^\varepsilon\|^2]+ L^2\mathbb{E}[|Y_s^\varepsilon|^2]  + \eta_s^2\big)\mathrm{d}s\\&\quad+ \int_t^T e^{\lambda s + \mu A_s}\nu_s^2\mathrm{d}A_s - 2 \int_t^T e^{\lambda s + \mu A_s} \langle Y_s^\varepsilon, Z_s^\varepsilon \mathrm{d}W_s \rangle.
		\end{aligned}
	\end{equation}
	and (10). By setting $t=0$ and taking the expectation in (12), we obtain
	\begin{equation}
		\mathbb{E}\int_0^T e^{\lambda s + \mu A_s} \big(|Y_s^\varepsilon|^2 (\mathrm{d}s + \mathrm{d}A_s) + \|Z_s^\varepsilon\|^2 \mathrm{d}s\big) \leq C \mathcal{M}(T). 
	\end{equation}
	Inequality (11) follows from combining (10) and (13).
	
	\noindent\textbf{Proposition 6}: Let Assumptions I-V be satisfied. Then there exists a positive constant \(C\) such that for any time \(t \in [0, T]\),
	\begin{equation}
		\begin{aligned}
			&\text{(a)} \quad \mathbb{E} \int_0^T e^{\lambda s + \mu A_s} \left( |\nabla \varphi_\varepsilon(Y_s^\varepsilon)|^2 \mathrm{d}s + |\nabla \psi_\varepsilon(Y_s^\varepsilon)|^2 \mathrm{d}A_s \right) \leq C \mathcal{M}(T), \\
			&\text{(b)} \quad \mathbb{E} \int_0^T e^{\lambda s + \mu A_s} \left( \varphi(J_\varepsilon (Y_s^\varepsilon)) \mathrm{d}s + \psi(\hat{J}_\varepsilon (Y_s^\varepsilon)) \mathrm{d}A_s \right) \leq C \mathcal{M}(T), \\
			&\text{(c)} \quad \mathbb{E} e^{\lambda t + \mu A_t} \left( |Y_t^\varepsilon - J_\varepsilon (Y_t^\varepsilon)|^2 + |Y_t^\varepsilon - \hat{J}_\varepsilon(Y_t^\varepsilon)|^2 \right) \leq \varepsilon  C \mathcal{M}(T), \\
			&\text{(d)} \quad \mathbb{E} e^{\lambda t + \mu A_t} \left( \varphi(J_\varepsilon (Y_t^\varepsilon)) + \psi(\hat{J}_\varepsilon (Y_t^\varepsilon)) \right) \leq C \mathcal{M}(T).
		\end{aligned}
	\end{equation}
	
	\textbf{Proof.} A key component of the proof is the stochastic subdifferential inequality developed by \cite{MR1642656}, which we use for our proof.
	
	We begin by presenting this subdifferential inequality.
	\begin{equation*}
		\begin{aligned}
			&e^{\lambda {t_{i+1} } + \mu A_{t_{i+1} }} \varphi_\varepsilon (Y_{t_{i+1} }^\varepsilon) 
			\geq \big(e^{\lambda t_{i+1} + \mu A_{t_{i+1} }} - e^{\lambda t_i + \mu A_{t_i } }\big) \varphi_\varepsilon(Y_{t_{i+1} }^\varepsilon) 
			+ e^{\lambda t_i  + \mu A_{t_i } }\varphi_\varepsilon(Y_{t_i } ^\varepsilon) \\
			&\quad + e^{\lambda t_i  + \mu A_{t_i } } \langle \nabla \varphi_\varepsilon(Y_{t_i } ^\varepsilon), Y_{t_{i+1} }^\varepsilon - Y_{t_i } ^\varepsilon \rangle,
		\end{aligned}
	\end{equation*}
	where \(t = t_0 < t_1 < t_2 < \cdots < T\) with \(t_{i+1} - t_i = \frac{1}{n}\). \\
	By summing up over \(i\) and passing to the limit as \(n \to \infty\), we get
	\[
	e^{\lambda T + \mu A_T} \varphi_\varepsilon(\xi) 
	\geq e^{\lambda t + \mu A_t} \varphi_\varepsilon(Y_t^\varepsilon) 
	+ \int_t^T e^{\lambda s + \mu A_s} \langle \nabla \varphi_\varepsilon(Y_s^\varepsilon), \mathrm{d}Y_s^\varepsilon \rangle +  \int_t^T \varphi_\varepsilon(Y_s^\varepsilon) \mathrm{d}(e^{\lambda s + \mu A_s}).
	\]
	Similarly, 
	\[
	e^{\lambda T + \mu A_T} \psi_\varepsilon(\xi) 
	\geq e^{\lambda t+ \mu A_t} \psi_\varepsilon(Y_t^\varepsilon) 
	+  \int_t^T e^{\lambda s + \mu A_s} \langle \nabla \psi_\varepsilon(Y_s^\varepsilon), \mathrm{d}Y_s^\varepsilon \rangle +  \int_t^T \psi_\varepsilon(Y_s^\varepsilon) \mathrm{d}(e^{\lambda s + \mu A_s}).
	\]
	Summing up the above and making use of (9), we get:
	\begin{equation}
		\begin{aligned}
			&e^{\lambda t + \mu A_t} \left( \varphi_\varepsilon(Y_t^\varepsilon) + \psi_\varepsilon(Y_t^\varepsilon) \right) 
			+  \int_t^T e^{\lambda s + \mu A_s} \left( |\nabla \varphi_\varepsilon(Y_s^\varepsilon)|^2 \mathrm{d}s + |\nabla \psi_\varepsilon(Y_s^\varepsilon)|^2 \mathrm{d}A_s \right) \\
			&\quad +  \int_t^T e^{\lambda s + \mu A_s} \left( \varphi_\varepsilon(Y_s^\varepsilon) + \psi_\varepsilon(Y_s^\varepsilon) \right) (\lambda \mathrm{d}s + \mu \mathrm{d}A_s)\\
			&\quad + \int_t^T e^{\lambda s + \mu A_s} \langle \nabla \varphi_\varepsilon (Y_s^\varepsilon) , \nabla \psi_\varepsilon (Y_s^\varepsilon)\rangle (\mathrm{d}s +\mathrm{d}A_s) \\
			&\leq e^{\lambda T + \mu A_T} \left( \varphi_\varepsilon(\xi) + \psi_\varepsilon(\xi) \right) 
			+ \int_t^T e^{\lambda s + \mu A_s} \langle \nabla \varphi_\varepsilon (Y_s^\varepsilon), F(s, Y_s^\varepsilon, Z_s^\varepsilon,\mathcal{L}(Y_s^\varepsilon, Z_s^\varepsilon)) \rangle \mathrm{d}s\\
			&\quad + \int_t^T e^{\lambda s + \mu A_s} \langle \nabla \varphi_\varepsilon (Y_s^\varepsilon), G(s, Y_s^\varepsilon) \rangle \mathrm{d}A_s + \int_t^T e^{\lambda s + \mu A_s} \langle \nabla \psi_\varepsilon (Y_s^\varepsilon), F(s, Y_s^\varepsilon, Z_s^\varepsilon,\mathcal{L}(Y_s^\varepsilon, Z_s^\varepsilon)) \rangle \mathrm{d}s \\
			&\quad +  \int_t^T e^{\lambda s + \mu A_s} \langle \nabla \psi_\varepsilon (Y_s^\varepsilon), G(s, Y_s^\varepsilon) \rangle \mathrm{d}A_s -  \int_t^T e^{\lambda s + \mu A_s} \langle \nabla \varphi_\varepsilon (Y_s^\varepsilon) + \nabla \psi_\varepsilon (Y_s^\varepsilon), Z_s^\varepsilon \mathrm{d}W_s \rangle .
		\end{aligned}
	\end{equation}
	With Assumption IV and the properties of the $W_2$-distance, we obtain
	\begin{equation*}
		\begin{aligned}
			&\frac{1}{2\varepsilon} |y - J_\varepsilon (y)|^2 \leq \varphi_\varepsilon (y), \quad \frac{1}{2\varepsilon} |y - \hat{J}_\varepsilon (y)|^2 \leq \psi_\varepsilon (y),\\
			&\varphi_\varepsilon (J_\varepsilon (y)) \leq \varphi_\varepsilon (y), \quad \psi_\varepsilon (\hat{J}_\varepsilon (y)) \leq \psi_\varepsilon (y),\\
			&\varphi_\varepsilon (\xi) \leq \varphi (\xi), \quad \psi_\varepsilon (\xi) \leq \psi (\xi),\\
			&|\langle \nabla \varphi_\varepsilon (Y_s), F(s, Y_s, Z_s,\mathcal{L}(Y_s,Z_s)) \rangle|\\&\leq \dfrac{1}{4}|\nabla \varphi_\varepsilon (Y_s)|^2 + 4L^2(|Y_s|^2+ \|Z_s\|^2+2\mathbb{E}[\|Z_s\|^2]+2\mathbb{E}[|Y_s|^2]) + 4\eta_s^2.\qquad\qquad\qquad\\
		\end{aligned}	
	\end{equation*}
	Considering Assumption V, we get:
	\begin{equation*}
		\begin{aligned}
			&\langle \nabla \psi_\varepsilon (Y_s), F(s, Y_s, Z_s, \mathcal{L}(Y_s,Z_s)) \rangle \leq
			\langle \nabla \varphi_\varepsilon (Y_s), F(s, Y_s, Z_s, \mathcal{L}(Y_s,Z_s)) \rangle^+ \\&\leq |\langle \nabla \varphi_\varepsilon (Y_s), F(s, Y_s, Z_s, \mu) \rangle|\\
			&\leq  \dfrac{1}{4}|\nabla \varphi_\varepsilon (Y_s)|^2 + 4L^2(|Y_s|^2+ \|Z_s\|^2+2\mathbb{E}[\|Z_s\|^2]+2\mathbb{E}[|Y_s|^2]) + 4\eta_s^2.
		\end{aligned}
	\end{equation*}
	Similarly,
	\begin{equation*}
		\begin{aligned}
			&\langle \nabla \varphi_\varepsilon (Y_s), G(s, Y_s) \rangle \leq \langle \nabla \psi_\varepsilon (Y_s), G(s, Y_s) \rangle^+ \leq |\langle \nabla \psi_\varepsilon (Y_s), G(s, Y_s) \rangle |\\&
			\leq \frac{1}{4} |\nabla \psi_\varepsilon (Y_s)|^2+ 2L^2 |Y_s|^2 + 2\nu_s^2. 
		\end{aligned}
	\end{equation*}
	Using the above inequality and employing a similar approach in Proposition 4, we obtain the conclusion of Proposition 5.
	
	\noindent\textbf{Proposition 7}:  Let Assumptions I-V be satisfied. Then
	\begin{equation}
		\begin{aligned}
			&\mathbb{E} \int_0^T e^{\lambda s + \mu A_s} \big( |Y_s^\varepsilon - Y_s^\delta|^2 (\mathrm{d}s + \mathrm{d}A_s) + \|Z_s^\varepsilon - Z_s^\delta\|^2 \mathrm{d}s \big) 
			+ \mathbb{E} \sup_{0 \leq s \leq T } e^{\lambda s + \mu A_s} |Y_s^\varepsilon - Y_s^\delta|^2 \\
			&\leq C (\varepsilon + \delta) \mathcal{M}(T). 
		\end{aligned}
	\end{equation}
	
	\textbf{Proof.} 	Similar to the previous derivation, we have 
	\begin{equation*}
		\begin{aligned}
			&2 \langle Y_s^\varepsilon - Y_s^\delta,  F(s, Y_s^\varepsilon, Z_s^\varepsilon, \mathcal{L}(Y_s^\varepsilon, Z_s^\varepsilon)) - F(s, Y_s^\delta, Z_s^\delta, \mathcal{L}(Y_s^\delta, Z_s^\delta) \rangle \\
			&\leq16L^2|Y_s^\varepsilon - Y_s^\delta|^2+ \frac{1}{4}|Z_s^\varepsilon - Z_s^\delta|^2+ \dfrac{1}{4}\mathbb{E}[\|Z_s^\varepsilon - Z_s^\delta\|^2]  + L^2\mathbb{E}[|Y_s^\varepsilon - Y_s^\delta|^2],\\
			&2 \langle Y_s^\varepsilon - Y_s^\delta, G(s, Y_s^\varepsilon) - G(s, Y_s^\delta) \rangle 
			\leq (K^2+1)|Y_s^\varepsilon - Y_s^\delta|^2.
		\end{aligned}
	\end{equation*}
	Given \(\nabla\varphi_\varepsilon=\dfrac{Y_s^\varepsilon-J_\varepsilon (Y_s^\varepsilon)}{\varepsilon}\) and using the property of operator \(\partial \varphi\) as discussed in \cite{MR274683}, which is $\nabla \varphi_\varepsilon (Y_t^\varepsilon) \in \partial \varphi (J_\varepsilon (Y_t^\varepsilon))$, we obtain
	\begin{equation*}
		\begin{aligned}
			&0 \leq \langle \nabla \varphi_\varepsilon (Y_s^\varepsilon) - \nabla \varphi_\delta (Y_s^\delta), J_\varepsilon (Y_s^\varepsilon) - J_\delta (Y_s^\delta) \rangle\\
			&= \langle \nabla \varphi_\varepsilon (Y_s^\varepsilon) - \nabla \varphi_\delta (Y_s^\delta), Y_s^\varepsilon - Y_s^\delta \rangle 
			- \varepsilon |\nabla \varphi_\varepsilon (Y_s^\varepsilon)|^2 - \delta |\nabla \varphi_\delta (Y_s^\delta)|^2\\&\quad + (\varepsilon + \delta) \langle \nabla \varphi_\varepsilon (Y_s^\varepsilon), \nabla \varphi_\delta (Y_s^\delta) \rangle,
		\end{aligned}
	\end{equation*}
	which leads to
	\begin{equation*}
		\langle \nabla \varphi_\varepsilon (Y_s^\varepsilon) - \nabla \varphi_\delta (Y_s^\delta), Y_s^\varepsilon - Y_s^\delta \rangle \geq - (\varepsilon + \delta) \langle \nabla \varphi_\varepsilon (Y_s^\varepsilon), \nabla \varphi_\delta (Y_s^\delta) \rangle	.
	\end{equation*}
	Similarly,
	\[
	\langle \nabla \psi_\varepsilon (Y_s^\varepsilon) - \nabla \psi_\delta (Y_s^\delta), Y_s^\varepsilon - Y_s^\delta \rangle 
	\geq - (\varepsilon + \delta) \langle \nabla \psi_\varepsilon (Y_s^\varepsilon), \nabla \psi_\delta (Y_s^\delta) \rangle.
	\]
	From Itô's formula and the above inequalities, we have
	\begin{equation*}
		\begin{aligned}
			&e^{\lambda t + \mu A_t} |Y_t^\varepsilon - Y_t^\delta|^2+ \int_t^T e^{\lambda s + \mu A_s} |Y_s^\varepsilon - Y_s^\delta|^2 (\lambda \mathrm{d}s + \mu \mathrm{d}A_s)+ \int_t^T e^{\lambda s + \mu A_s} \|Z_s^\varepsilon - Z_s^\delta\|^2 \mathrm{d}s \\
			&\quad- 2 \int_t^T e^{\lambda s + \mu A_s}  (\varepsilon + \delta) \langle \nabla \varphi_\varepsilon (Y_s^\varepsilon), \nabla \varphi_\delta (Y_s^\delta) \rangle \mathrm{d}s\\
			&\quad - 2 \int_t^T e^{\lambda s + \mu A_s} (\varepsilon + \delta) \langle \nabla \psi_\varepsilon (Y_s^\varepsilon), \nabla \psi_\delta (Y_s^\delta) \rangle \mathrm{d}A_s \\
			&\leq \int_t^T e^{\lambda s + \mu A_s}  \big((2+9L^2)|Y_s^\varepsilon - Y_s^\delta|^2+ \frac{1}{4}\|Z_s^\varepsilon - Z_s^\delta\|^2\big) \mathrm{d}s \\
			&\quad + \int_t^T e^{\lambda s + \mu A_s}\big(\dfrac{1}{4}\mathbb{E}[\|Z_s^\varepsilon - Z_s^\delta\|^2]  + L^2\mathbb{E}[|Y_s^\varepsilon - Y_s^\delta|^2]\big)\mathrm{d}s \\
			&\quad + \int_t^T e^{\lambda s + \mu A_s}  (K^2+1)|Y_s^\varepsilon - Y_s^\delta|^2 \mathrm{d}A_s - 2 \int_t^T e^{\lambda s + \mu A_s} \langle Y_s^\varepsilon - Y_s^\delta, (Z_s^\varepsilon - Z_s^\delta) \mathrm{d}W_s \rangle.\\
		\end{aligned}
	\end{equation*}
	From (14)-(a), 
	\begin{equation*}
		\begin{aligned}
			&2(\varepsilon + \delta)\mathbb{E} \int_t^T e^{\lambda s + \mu A_s} [\langle \nabla \varphi_\varepsilon (Y_s^\varepsilon), \nabla \varphi_\delta (Y_s^\delta) \rangle \mathrm{d}s + \langle \nabla \psi_\varepsilon (Y_s^\varepsilon), \nabla \psi_\delta (Y_s^\delta) \rangle \mathrm{d}A_s]\\&\leq C(\varepsilon + \delta) \mathcal{M}(T).
		\end{aligned}
	\end{equation*}
	Following an argument similar to that in Proposition 4, we obtain (15).
	
	\subsubsection{Proof of Theorem 3}
	
	From Proposition 5, we have
	\[
	\begin{aligned}
		&\exists Y \in \mathcal{S}_k^{ \lambda, \mu} \cap H_k^{\lambda, \mu} \cap \tilde{\mathcal{H}}_k^{\lambda, \mu} , \quad \exists Z \in \mathcal{H}_{k \times d}^{ \lambda, \mu}, \\
		&\lim_{\varepsilon  \searrow  0} Y^\varepsilon = Y\quad \text{in } \mathcal{S}_k^{ \lambda, \mu} \cap \mathcal{H}_k^{\lambda, \mu} \cap \tilde{\mathcal{H}}_k^{\lambda, \mu}, \\
		&\lim_{\varepsilon  \searrow  0} Z^\varepsilon = Z\quad \text{in }  \mathcal{H}_{k \times d}^{ \lambda, \mu}.
	\end{aligned}
	\]
	From (14)-(a) and (14)-(c), for any time \(t \in [ 0, T]\), it follows that
	\[
	\lim_{\varepsilon  \searrow  0} J_\varepsilon (Y^\varepsilon) = Y \quad \text{in } \mathcal{H}_k^{\lambda, \mu},
	\quad \lim_{\varepsilon  \searrow  0} \hat{J}_\varepsilon (Y^\varepsilon) = Y \quad \text{in } \tilde{\mathcal{H}}_k^{\lambda, \mu},
	\]
	\[
	\lim_{\varepsilon  \searrow 0} \mathbb{E} e^{\lambda t + \mu A_t} |J_\varepsilon (Y_t^\varepsilon) - Y_t|^2 = 0,
	\quad \lim_{\varepsilon  \searrow  0} \mathbb{E} e^{\lambda t + \mu A_t} |\hat{J}_\varepsilon (Y_t^\varepsilon) - Y_t|^2 = 0.
	\]
	Using Fatou’s Lemma, (14)-(b), (14)-(d) and and the lower semicontinuity of $\varphi\, \psi$, we establish (c) of Definition 1.
	
	We define $U^\varepsilon := \nabla \varphi_\varepsilon (Y^\varepsilon)$ and $V^\varepsilon := \nabla \psi_\varepsilon (Y^\varepsilon)$. From (14)-(a), it follows that
	\[
	\mathbb{E} \int_t^T e^{\lambda s + \mu A_s} \big(|U^\varepsilon|^2 \mathrm{d}s + |V^\varepsilon|^2 \mathrm{d}A_s \big) \leq C M(T).
	\]
	Hence there exist \(U \in \mathcal{H}_k^{\lambda, \mu}\) and \(V \in \tilde{\mathcal{H}}_k^{\lambda, \mu}\) such that for a subsequence \(\varepsilon_n \searrow   0\),
	\[
	U^{\varepsilon_n} \rightharpoonup U, \quad \text{weakly in Hilbert space } \mathcal{H}_k^{\lambda, \mu},
	\]
	\[
	V^{\varepsilon_n} \rightharpoonup V, \quad \text{weakly in Hilbert space } \tilde{\mathcal{H}}_k^{\lambda, \mu}.
	\]
	Then
	\[
	\mathbb{E} \int_t^T e^{\lambda s + \mu A_s} \big(|U|^2 \mathrm{d}s + |V|^2 \mathrm{d}A_s \big)
	\leq \liminf_{n \to \infty} \mathbb{E} \int_t^T e^{\lambda s + \mu A_s} \big(|U^{\varepsilon_n}|^2 \mathrm{d}s + |V^{\varepsilon_n}|^2 \mathrm{d}A_s \big)
	\leq C \mathcal{M} (t,T),
	\]
	where $\mathcal{M}(t, T)$ shares the same structure as $\mathcal{M}(T)$, with the only difference being that the lower bound of the integral is changed to $t$. Taking the limit on both sides of equation (9), we obtain (e) of Definition 1.
	
	Let \(u \in  \mathcal{H}_k^{\lambda, \mu}, v \in \tilde{ \mathcal{H}}_k^{\lambda, \mu}\). \(\nabla \varphi_\varepsilon (Y_t^\varepsilon) \in \partial \varphi (J_\varepsilon (Y_t^\varepsilon))\) and \(\nabla \psi_\varepsilon (Y_t^\varepsilon) \in \partial \psi (\hat{J}_\varepsilon (Y_t^\varepsilon))\), \(\forall \, t \in [0,T]\). Then as signed measures on \(\Omega \times [0, T]\),
	\begin{equation*}
		\begin{aligned}
			&e^{\lambda s + \mu A_s} \langle U_s^\varepsilon, u_s - J_\varepsilon (Y_s^\varepsilon)\rangle \mathbb{P} (\mathrm{d}\omega) \otimes \mathrm{d}s + e^{\lambda s + \mu A_s} \varphi \big(J_\varepsilon (Y_s^\varepsilon)\big) \, \mathbb{P} (\mathrm{d}\omega) \otimes \mathrm{d}s\\
			&\leq e^{\lambda s + \mu A_s} \varphi (u_s) \, \mathbb{P} (\mathrm{d}\omega) \otimes \mathrm{d}s,\\
			&e^{\lambda s + \mu A_s} \langle V_s^\varepsilon, v_s - \hat{J}_\varepsilon (Y_s^\varepsilon)\rangle \mathbb{P} (\mathrm{d}\omega) \otimes A (\omega, \mathrm{d}s)+ e^{\lambda s + \mu A_s} \psi \big(\hat{J}_\varepsilon (Y_s^\varepsilon)\big) \, \mathbb{P} (\mathrm{d}\omega) \otimes A (\omega, \mathrm{d}s)\\
			&\leq e^{\lambda s + \mu A_s} \psi (v_s) \, \mathbb{P} (\mathrm{d}\omega) \otimes A (\omega, \mathrm{d}s).
		\end{aligned}
	\end{equation*}
	Taking the limit in the two inequalities above, we obtain (d) of Definition 1, thereby completing the proof of existence.
	
	By Proposition 2, if the difference between the two terminal values is zero, then the right-hand side of inequality (3) vanishes, which implies the uniqueness of the solution.
	
	\section{BSVIs with locally Lipschitz and non-linear growth conditions}
	We reformulate (iii)-(iv) in Assumption I as (iii')-(iv'), respectively, and denote the resulting modification as Assumption I'. 
	\begin{enumerate}[label=(\roman*)]
		\item[(iii')] $\forall N \in \mathcal{N}, \exists r_N \in \mathbb{R}$ such that, for any $y,y',z,z'$ satisfying $|y|,|y'|,\|z\|,\|z'\|\leq N$ and for any $\Pi,\Pi'$ satisfying $W_2(\Pi,\delta_0),W_2(\Pi',\delta_0)\leq N$,
		$$|F(t, y,z,\Pi) - F(t,y',z',\Pi')| \leq r_N( |y-y'|+|z - z'|+W_2(\Pi,\Pi')) ,$$
		\item[(iv')] for $\alpha \in [0,1)$, there exists $L\geq 0$ such that \(|F(t, y, z,\Pi)| \leq L(|y|^\alpha+|z|^\alpha+W_2(\Pi,\delta_0))+\eta_t \).
	\end{enumerate}	
	
	We present a technical lemma whose proof similar to the arguments used in Lemma 3.2 of \cite{MR4637498}. 
	
	\noindent\textbf{Lemma 8}: 
	Under Assumption I', there exists a sequence $\{F_n\}_{n=1}^{\infty}$ such that
	\begin{itemize}
		\item[(i)] $\sup_{n\geq 0}|F_n(t,y,z,\Pi)|\leq|F(t,y,z,\Pi)|\leq L(|y|^\alpha+|z|^\alpha+W_2(\Pi,\delta_0))+\eta_t,$
		\item[(ii)] $\forall n, F_n$ is globally Lipschitz, e.g
		$$|F_n(t, y,z,\Pi) - F_n(t,y',z',\Pi')| \leq L_n( |y-y'|+|z - z'|+W_2(\Pi,\Pi')) ,$$
		\item[(iii)] $\forall N \in \mathcal{N},\; \rho_N(F_n - F) \to 0$ as $n \to \infty$, where $$ \rho_N^2(f) := \mathbb{E} \left[ \int_0^T \sup_{|y|, |z|,W_2(\Pi,\delta_0)  \leq N} e^{\lambda t + \mu A_t} |f(t, y, z,\Pi)|^2 \, \mathrm{d}t \right].$$
	\end{itemize}
	
	\text{Proof. }Let $\phi^1_n : \mathbb{R} \to \mathbb{R}_+$ be a sequence of smooth functions such that $0 \leq \phi^1_n \leq 1$, $\phi^1_n(x) = 1$ for $|x| \leq n$, and $\phi^1_n(x) = 0$ for $|x| \geq n + 1$. Similarly, we define $\phi^2_n : \mathbb{R}^d \to \mathbb{R}_+$ and $\phi^3_n : \mathbb{R}^{d \times k} \to \mathbb{R}_+$. We define $F_n(t, y, z,\Pi) := F(t, y, z, \Pi)\phi^2_n(y)\phi^3_n(z)\psi^1_n(W_2(\Pi,\delta_0))$. $F_n(t, y, z, \Pi)$ satisfies the claims.
	
	\noindent\textbf{Theorem 9}: Let Assumptions I' and II–V be satisfied. Then there exists a unique solution to BSVI (1).
	
	\textbf{Proof.} Let $F_n$ be the approximation sequence of $F$ in Lemma 8. In view of (ii) in Lemma 8, we set $\epsilon=\dfrac{1}{n}$ in (9) and apply an argument similar to that in section 2.2 to show that the following equation admits a unique solution 
	\begin{equation}
		\begin{aligned}
			&Y_t^{\frac{1}{n}}+ \int_{ t}^T \nabla \varphi_{\frac{1}{n}}(Y_s^{\frac{1}{n}}) \mathrm{d}s + \int_{t}^T \nabla \psi_{\frac{1}{n}}(Y_s^{\frac{1}{n}}) \mathrm{d}A_s 
			= \xi + \int_{ t}^T F_n(s, Y_s^{\frac{1}{n}}, Z_s^{\frac{1}{n}},\mathcal{L}(Y_s^{\frac{1}{n}}, Z_s^{\frac{1}{n}})) \mathrm{d}s \\
			&\quad+ \int_{ t}^T G(s, Y_s^{\frac{1}{n}}) \mathrm{d}A_s 
			- \int_{ t}^T Z_s^{\frac{1}{n}} \mathrm{d}W_s, \quad \forall t \in[0,T].
		\end{aligned}
	\end{equation}
	
	Step 1): We first proof the following estimation,
	\begin{equation}
		\begin{aligned}
			&\mathbb{E} \left[ \sup_{0 \leq s \leq T} e^{\lambda s + \mu A_s} |Y_s^{\frac{1}{n}}|^2 
			+ \int_0^T e^{\lambda s + \mu A_s} \big(|Y_s^{\frac{1}{n}}|^2 + \|Z_s^{\frac{1}{n}}\|^2\big) \mathrm{d}s
			+ \int_0^T e^{\lambda s + \mu A_s} |Y_s^{\frac{1}{n}}|^2 \mathrm{d}A_s \right] \leq C \mathcal{M}(T).
		\end{aligned}
	\end{equation}
	
    First, we present the proof of existence. Analogous to the proof of Proposition 4, we have
	\begin{equation*}
			2 \langle Y_s^{\frac{1}{n}}, G(s, Y_s^{\frac{1}{n}}) \rangle \leq(2+K^2)|Y_s^{\frac{1}{n}}|^2 + \nu_s^2.\qquad\qquad\qquad\qquad\qquad
	\end{equation*}
    Considering the estimates of the generator $F$ in its integral form, we have
    \begin{equation*}
		\begin{aligned}
			&2 \int_{ t}^Te^{\lambda s + \mu A_s}\langle Y_s^{\frac{1}{n}}, F_n(s, Y_s^{\frac{1}{n}}, Z_s^{\frac{1}{n}},\mathcal{L}(Y_s^{\frac{1}{n}}, Z_s^{\frac{1}{n}})) \rangle \mathrm{d}s\\
            &\leq \int_{t}^Te^{\lambda s + \mu A_s}\big[p|Y_s^{\frac{1}{n}}|^2+\frac{1}{p}|F_n(s, Y_s^{\frac{1}{n}}, Z_s^{\frac{1}{n}},\mathcal{L}(Y_s^{\frac{1}{n}}, Z_s^{\frac{1}{n}}))|^2\big]\mathrm{d}s\\
            & \leq \int_{t}^Te^{\lambda s + \mu A_s}\big[p|Y_s^{\frac{1}{n}}|^2+\frac{1}{p}|L(|Y_s^{\frac{1}{n}}|^\alpha+\|Z_s^{\frac{1}{n}}\|^\alpha+W_2(\Pi,\delta_0))+\eta_t|^2\big]\mathrm{d}s\\
            & \leq \int_{t}^Te^{\lambda s + \mu A_s}\big[p|Y_s^{\frac{1}{n}}|^2+\frac{1}{p}\big(qL^2(|Y_s^{\frac{1}{n}}|^\alpha+\|Z_s^{\frac{1}{n}}\|^\alpha+W_2(\Pi,\delta_0))^2+\frac{1}{q}\eta_t^2\big)\big]\mathrm{d}s\\
            & \leq \int_{t}^Te^{\lambda s + \mu A_s}\big[p|Y_s^{\frac{1}{n}}|^2+\frac{1}{pq}\eta_t^2+\frac{3L^2q}{p}(|Y_s^{\frac{1}{n}}|^{2\alpha}+\|Z_s^{\frac{1}{n}}\|^{2\alpha})+\frac{6L^2q}{p}(\mathbb{E}|Y_s^{\frac{1}{n}}|^2+  \mathbb{E}\|Z_s^{\frac{1}{n}}\|^2)\big]\mathrm{d}s.
		\end{aligned}
	\end{equation*}
    Applying the integral form of Hölder's inequality to the terms involving $|Y_s^{\frac{1}{n}}|^{2\alpha}$, we obtain
    \begin{equation*}
        \begin{aligned}
            &\int_{t}^Te^{\lambda s + \mu A_s}|Y_s^{\frac{1}{n}}|^{2\alpha}\mathrm{d}s\\
            &\leq\big[\int_{t}^T\big(e^{\alpha(\lambda s + \mu A_s)}|Y_s^{\frac{1}{n}}|^{2\alpha}\big)^{\frac{1}{\alpha}}\mathrm{d}s\big]^\alpha\big[\int_{t}^T\big(e^{(1-\alpha)(\lambda s + \mu A_s)}\big)^{\frac{1}{1-\alpha}}\mathrm{d}s\big]^{1-\alpha}\\
            &\leq\big[\int_{t}^Te^{\lambda s + \mu A_s}|Y_s^{\frac{1}{n}}|^{2}\mathrm{d}s\big]^\alpha\big(\int_{t}^Te^{\lambda s + \mu A_s}\mathrm{d}s\big)^{1-\alpha}\\
            &\leq Te^{\lambda T + \mu A_T}\int_{t}^Te^{\lambda s + \mu A_s}|Y_s^{\frac{1}{n}}|^{2}\mathrm{d}s.
        \end{aligned}
    \end{equation*}
    Following a similar procedure for the terms containing $|Z_s^{\frac{1}{n}}|^{2\alpha}$ and setting $q=\frac{1}{Te^{\lambda T + \mu A_T}}, p=12L^2$, the previous expressions can be consolidated as follows:
    \begin{equation*}
		\begin{aligned}
			&2 \int_{ t}^Te^{\lambda s + \mu A_s}\langle Y_s^{\frac{1}{n}}, F_n(s, Y_s^{\frac{1}{n}}, Z_s^{\frac{1}{n}},\mathcal{L}(Y_s^{\frac{1}{n}}, Z_s^{\frac{1}{n}})) \rangle \mathrm{d}s\\
            & \leq \int_{ t}^T e^{\lambda s + \mu A_s} \big[(\frac{1}{4}+12L^2)|Y_s^{\frac{1}{n}}|^2+\frac{1}{4}\|Z_s^{\frac{1}{n}}\|^2+\frac{1}{2} \mathbb{E}|Y_s^{\frac{1}{n}}|^2+\frac{1}{2}  \mathbb{E}\|Z_s^{\frac{1}{n}}\|^2+\frac{12L^2}{ T}\eta_s^2\big]\mathrm{d}s.
		\end{aligned}
	\end{equation*}
    Using the Ito’s formula and the above inequality, we get
	\begin{equation*}
		\begin{aligned}
			&e^{\lambda t + \mu A_{t }} |Y_{t}^{\frac{1}{n}}|^2 
			+ \int_{ t}^T e^{\lambda s + \mu A_s} \big(|Y_s^{\frac{1}{n}}|^2 (\lambda \mathrm{d}s + \mu \mathrm{d}A_s) 
			+ \|Z_s^{\frac{1}{n}}\|^2 \mathrm{d}s\big) \\
			&\quad + 2 \int_{ t}^T e^{\lambda s + \mu A_s} \big(\langle Y_s^{\frac{1}{n}}, \nabla \varphi_{\frac{1}{n}}(Y_s^{\frac{1}{n}}) \rangle \mathrm{d}s
			+ \langle Y_s^{\frac{1}{n}}, \nabla \varphi_{\frac{1}{n}}(Y_s^{\frac{1}{n}}) \rangle \mu \mathrm{d}A_s \big) \\
			&\leq e^{\lambda T + \mu A_T} |\xi|^2 + \int_{ t}^T e^{\lambda s + \mu A_s} \big((2+K^2)|Y_s^{\frac{1}{n}}|^2 + \nu_s^2\big)\mathrm{d}A_s - 2 \int_{ t}^T e^{\lambda s + \mu A_s} \langle Y_s^{\frac{1}{n}}, Z_s^{\frac{1}{n}} \mathrm{d}W_s \rangle \\
            & \leq \int_{ t}^T e^{\lambda s + \mu A_s} \big[(\frac{1}{4}+12L^2)|Y_s^{\frac{1}{n}}|^2+\frac{1}{4}\|Z_s^{\frac{1}{n}}\|^2+\frac{1}{2} \mathbb{E}|Y_s^{\frac{1}{n}}|^2+\frac{1}{2}  \mathbb{E}\|Z_s^{\frac{1}{n}}\|^2+\frac{12L^2}{ T}\eta_s^2\big]\mathrm{d}s.
		\end{aligned}
	\end{equation*}
	By applying the method used to obtain (10) and (13) and setting $t = 0$, it follows that
	\begin{equation*}
        \begin{aligned}
            &\mathbb{E}\int_0^T e^{\lambda s + \mu A_s} \big(|Y_s^{\frac{1}{n}}|^2 (\mathrm{d}s + \mathrm{d}A_s) + \|Z_s^{\frac{1}{n}}\|^2 \mathrm{d}s\big) \leq C \mathcal{M}(T) , \\
            &\mathbb{E} \sup_{0 \leq s \leq T} e^{\lambda s + \mu A_s} | Y_s^{\frac{1}{n}} |^2 \leq C \mathcal{M}(T) .
        \end{aligned}
	\end{equation*}
	Consequently, we obtain (18). 
	
	Step 2): An argument analogous to that of Proposition 6, combined with the fact that $F_n$ is globally Lipschitz, yields the following result. We omit the proof here. 
	
	There exists a positive constant \(C\) such that for any time \(t \in [0, T]\),
	\begin{equation*}
		\begin{aligned}
			&\mathbb{E} \int_0^T e^{\lambda s + \mu A_s} \left( |\nabla \varphi_{\frac{1}{n}}(Y_s^{\frac{1}{n}})|^2 \mathrm{d}s + |\nabla \psi_{\frac{1}{n}}(Y_s^{\frac{1}{n}})|^2 \mathrm{d}A_s \right) \leq C_n \mathcal{M}(T), \\
			&\mathbb{E} \int_0^T e^{\lambda s + \mu A_s} \left( \varphi(J_{\frac{1}{n}} (Y_s^{\frac{1}{n}})) \mathrm{d}s + \psi(\hat{J}_{\frac{1}{n}} (Y_s^{\frac{1}{n}})) \mathrm{d}A_s \right) \leq C_n \mathcal{M}(T), \\
			&\mathbb{E} e^{\lambda t + \mu A_t} \left( |Y_t^{\frac{1}{n}} - J_{\frac{1}{n}} (Y_t^{\frac{1}{n}})|^2 + |Y_t^{\frac{1}{n}} - \hat{J}_{\frac{1}{n}}(Y_t^{\frac{1}{n}})|^2 \right) \leq {\frac{C_n}{n}}\mathcal{M}(T), \\
			&\mathbb{E} e^{\lambda t + \mu A_t} \left( \varphi(J_{\frac{1}{n}} (Y_t^{\frac{1}{n}})) + \psi(\hat{J}_{\frac{1}{n}} (Y_t^{\frac{1}{n}})) \right) \leq C_n \mathcal{M}(T).
		\end{aligned}
	\end{equation*}
	
	Step 3): For any given $N$, we set $D_{n,m}^N:=\{(s,w): |Y_s^{\frac{1}{n}}|^2+|Y_s^{\frac{1}{m}}|^2+|Z_s^{\frac{1}{n}}|^2+|Z_s^{\frac{1}{m}}|^2\geq N^2\}$ and denote $\bar{D}_{n,m}^N=\Omega-D_{n,m}^N$. 
	
    By using Itô formula and taking intergal, we have
    \begin{equation*}
		\begin{aligned}
			&e^{\lambda t +\mu A_{t }} |Y_{t }^{\frac{1}{n}} - Y_{t }^{\frac{1}{m}}|^2 
			+ \int_{ t}^T e^{\lambda s + \mu A_s} |Y_s^{\frac{1}{n}} - Y_s^{\frac{1}{m}}|^2 (\lambda \mathrm{d}s + \mu \mathrm{d}A_s) +\int_{ t}^T e^{\lambda s + \mu A_s} \|Z_s^{\frac{1}{n}} - Z_s^{\frac{1}{m}}\|^2 \mathrm{d}s  \\
			&\quad + 2 \int_{ t}^T e^{\lambda s + \mu A_s} \langle Y_s^{\frac{1}{n}} - Y_s^{\frac{1}{m}}, \nabla \varphi_n (Y_s^{\frac{1}{n}}) - \nabla \varphi_m (Y_s^{\frac{1}{m}}) \rangle \mathrm{d}s \\&\quad + 2 \int_{ t}^T e^{\lambda s + \mu A_s} \langle Y_s^{\frac{1}{n}} - Y_s^{\frac{1}{m}}, \nabla \psi_n (Y_s^{\frac{1}{n}}) - \nabla \psi_m (Y_s^{\frac{1}{m}}) \rangle \mathrm{d}A_s \\
			&= 2 \int_{ t}^T e^{\lambda s + \mu A_s} \langle Y_s^{\frac{1}{n}} - Y_s^{\frac{1}{m}}, F_n(s, Y_s^{\frac{1}{n}}, Z_s^{\frac{1}{n}}, \mathcal{L}(Y_s^{\frac{1}{n}}, Z_s^{\frac{1}{n}})) - F_m(s, Y_s^{\frac{1}{m}}, Z_s^{\frac{1}{m}}, \mathcal{L}(Y_s^{\frac{1}{m}}, Z_s^{\frac{1}{m}})) \rangle \mathrm{d}s \\
			&\quad + 2 \int_{ t}^T e^{\lambda s + \mu A_s} \langle Y_s^{\frac{1}{n}} - Y_s^{\frac{1}{m}}, G(s, Y_s^{\frac{1}{n}}) - G(s, Y_s^{\frac{1}{m}}) \rangle \mathrm{d}A_s\\
            &\quad -2 \int_{ t}^T e^{\lambda s + \mu A_s} \langle \langle Y_s^{\frac{1}{n}} - Y_s^{\frac{1}{m}},(Z_s^{\frac{1}{n}} - Z_s^{\frac{1}{m}})dW_s\rangle\\
            &=:\Gamma_1+\Gamma_2+\Gamma_3+\Gamma_4+\Gamma_5+\Gamma_6,\\
		\end{aligned}
	\end{equation*}
	where,
    \begin{equation*}
		\begin{aligned}
			&\Gamma_1:= 2 \int_{ t}^T e^{\lambda s + \mu A_s} \langle Y_s^{\frac{1}{n}} - Y_s^{\frac{1}{m}},\Theta_1\rangle 1_{D_{n,m}^N} \mathrm{d}s ,\quad \Gamma_2:=2 \int_{ t}^T e^{\lambda s + \mu A_s} \langle Y_s^{\frac{1}{n}} - Y_s^{\frac{1}{m}},\Theta_2\rangle 1_{\bar{D}_{n,m}^N}\mathrm{d}s, \\
			&\Gamma_3:=2 \int_{ t}^T e^{\lambda s + \mu A_s} \langle Y_s^{\frac{1}{n}} - Y_s^{\frac{1}{m}},\Theta_3\rangle1_{\bar{D}_{n,m}^N}  \mathrm{d}s, \quad \Gamma_4:=2 \int_{ t}^T e^{\lambda s + \mu A_s} \langle Y_s^{\frac{1}{n}} - Y_s^{\frac{1}{m}},\Theta_4 \rangle1_{\bar{D}_{n,m}^N} \mathrm{d}s, \\
			&\Gamma5:=2 \int_{ t}^T e^{\lambda s + \mu A_s} \langle Y_s^{\frac{1}{n}} - Y_s^{\frac{1}{m}}, G(s, Y_s^{\frac{1}{n}}) - G(s, Y_s^{\frac{1}{m}}) \rangle \mathrm{d}A_s,\\
            &\Gamma_6=-2 \int_{ t}^T e^{\lambda s + \mu A_s} \langle \langle Y_s^{\frac{1}{n}} - Y_s^{\frac{1}{m}},(Z_s^{\frac{1}{n}} - Z_s^{\frac{1}{m}})dW_s\rangle,\\
            &\Theta_1:=F_n(s, Y_s^{\frac{1}{n}}, Z_s^{\frac{1}{n}}, \mathcal{L}(Y_s^{\frac{1}{n}}, Z_s^{\frac{1}{n}})) - F_m(s, Y_s^{\frac{1}{m}}, Z_s^{\frac{1}{m}}, \mathcal{L}(Y_s^{\frac{1}{m}}, Z_s^{\frac{1}{m}})),\\
            &\Theta_2:=F_n(s, Y_s^{\frac{1}{n}}, Z_s^{\frac{1}{n}}, \mathcal{L}(Y_s^{\frac{1}{n}}, Z_s^{\frac{1}{n}})) - F(s, Y_s^{\frac{1}{n}}, Z_s^{\frac{1}{n}}, \mathcal{L}(Y_s^{\frac{1}{n}}, Z_s^{\frac{1}{n}})) ,\\
            &\Theta_3:=F(s, Y_s^{\frac{1}{n}}, Z_s^{\frac{1}{n}}, \mathcal{L}(Y_s^{\frac{1}{n}}, Z_s^{\frac{1}{n}})) - F(s, Y_s^{\frac{1}{m}}, Z_s^{\frac{1}{m}}, \mathcal{L}(Y_s^{\frac{1}{m}}, Z_s^{\frac{1}{m}})),\\
            &\Theta_4:=F(s, Y_s^{\frac{1}{m}}, Z_s^{\frac{1}{m}}, \mathcal{L}(Y_s^{\frac{1}{m}}, Z_s^{\frac{1}{m}})) - F_m(s, Y_s^{\frac{1}{m}}, Z_s^{\frac{1}{m}}, \mathcal{L}(Y_s^{\frac{1}{m}}, Z_s^{\frac{1}{m}})).
		\end{aligned}
	\end{equation*}

    Utilizing the assumptions, properties of $D_{n,m}^N$, (18) and ${x}^\alpha\leq1+|x|, \alpha\in [0,1)$, we derive the following inequalities
    \begin{equation*}
		\begin{aligned}
			&\Gamma_1\leq \int_{ t}^T e^{\lambda s + \mu A_s} |Y_s^{\frac{1}{n}} - Y_s^{\frac{1}{m}}|^2\mathrm{d}s+20L^2\int_{ t}^T e^{\lambda s + \mu A_s} \eta_s^21_{D_{n,m}^N}\mathrm{d}s\\
			&\quad+10L^2\int_{ t}^T e^{\lambda s + \mu A_s}\mathbb{E}\big(|Y_s^{\frac{1}{n}}|^2+|Y_s^{\frac{1}{m}}|^2+\|Z_s^{\frac{1}{n}}\|^2+\|Z_s^{\frac{1}{m}}\|^2\big)1_{D_{n,m}^N}\mathrm{d}s\\
			&\quad+10L^2\int_{ t}^T e^{\lambda s + \mu A_s}\big(|Y_s^{\frac{1}{n}}|^{2\alpha}+|Y_s^{\frac{1}{m}}|^{2\alpha}+\|Z_s^{\frac{1}{n}}\|^{2\alpha}+\|Z_s^{\frac{1}{m}}\|^{2\alpha}\big)1_{D_{n,m}^N}\mathrm{d}s\\
            & \leq \int_{ t}^T e^{\lambda s + \mu A_s} |Y_s^{\frac{1}{n}} - Y_s^{\frac{1}{m}}|^2\mathrm{d}s\\
            &\quad+20L^2\int_{ t}^T e^{\lambda s + \mu A_s} \eta_s^2\dfrac{|Y_s^{\frac{1}{n}}|^2+|Y_s^{\frac{1}{m}}|^2+|Z_s^{\frac{1}{n}}|^2+|Z_s^{\frac{1}{m}}|^2}{N^2}1_{D_{n,m}^N}\mathrm{d}s\\
            &\quad+10L^2\int_{ t}^T \mathbb{E}\big(\sup_{0 \leq s \leq T} e^{\lambda s + \mu A_s}|Y_s^{\frac{1}{n}}|^2+\sup_{0 \leq s \leq T} e^{\lambda s + \mu A_s}|Y_s^{\frac{1}{m}}|^2+\int_{0}^Te^{\lambda r + \mu A_r}\|Z_r^{\frac{1}{n}}\|^2\mathrm{d}r\\&\qquad+\int_{0}^Te^{\lambda r + \mu A_r}\|Z_s^{\frac{1}{m}}\|^2\mathrm{d}r\big)\dfrac{|Y_s^{\frac{1}{n}}|^2+|Y_s^{\frac{1}{m}}|^2+|Z_s^{\frac{1}{n}}|^2+|Z_s^{\frac{1}{m}}|^2}{N^2}1_{D_{n,m}^N}\mathrm{d}s\\
            &\quad+10L^2\int_{ t}^T e^{\lambda s + \mu A_s}\big(4+|Y_s^{\frac{1}{n}}|^{2}+|Y_s^{\frac{1}{m}}|^{2}+\|Z_s^{\frac{1}{n}}\|^{2}+\|Z_s^{\frac{1}{m}}\|^{2}\big)\dfrac{|Y_s^{\frac{1}{n}}|^2+|Y_s^{\frac{1}{m}}|^2+|Z_s^{\frac{1}{n}}|^2+|Z_s^{\frac{1}{m}}|^2}{N^2}1_{D_{n,m}^N}\mathrm{d}s\\
            & \leq \int_{ t}^T e^{\lambda s + \mu A_s} |Y_s^{\frac{1}{n}} - Y_s^{\frac{1}{m}}|^2\mathrm{d}s\\
            &\quad+20L^2\big[\int_{ t}^T e^{\lambda s + \mu A_s} \eta_s^4\mathrm{d}s \big]^{\frac{1}{2}}\big[\int_{ t}^T e^{\lambda s + \mu A_s}(\dfrac{|Y_s^{\frac{1}{n}}|^2+|Y_s^{\frac{1}{m}}|^2+|Z_s^{\frac{1}{n}}|^2+|Z_s^{\frac{1}{m}}|^2}{N^2})^21_{D_{n,m}^N}\mathrm{d}s\big]^{\frac{1}{2}}\\
            &\quad +C\mathcal{M}(T)\mathbb{E}\int_{ t}^T e^{\lambda s + \mu A_s}\dfrac{|Y_s^{\frac{1}{n}}|^2+|Y_s^{\frac{1}{m}}|^2+|Z_s^{\frac{1}{n}}|^2+|Z_s^{\frac{1}{m}}|^2}{N^2}1_{D_{n,m}^N}\mathrm{d}s\\
            &\quad +\frac{40L^2}{N^2}\mathbb{E}\int_{ t}^T e^{\lambda s + \mu A_s}|Y_s^{\frac{1}{n}}|^{2}+|Y_s^{\frac{1}{m}}|^{2}+\|Z_s^{\frac{1}{n}}\|^{2}+\|Z_s^{\frac{1}{m}}\|^{2}1_{D_{n,m}^N}\mathrm{d}s\\
            &\quad +\frac{10L^2}{N^2}\mathbb{E}\int_{ t}^T e^{\lambda s + \mu A_s}\big(|Y_s^{\frac{1}{n}}|^{2}+|Y_s^{\frac{1}{m}}|^{2}+\|Z_s^{\frac{1}{n}}\|^{2}+\|Z_s^{\frac{1}{m}}\|^{2}\big)^21_{D_{n,m}^N}\mathrm{d}s.\\
            \end{aligned}
        \end{equation*}

    Applying Itô's formula to $e^{\lambda s + \mu A_s}|Y_s|^4$, and employing Hölder's inequality along with similar techniques as in the proof of (18), we can obtain results with a similar form but different orders. Since the derivation is rather lengthy, we omit the details here. Denoting the corresponding upper bound on the right-hand side as $\tilde{\mathcal{M}}(T)$, we have
    \begin{equation*}
        \mathbb{E}[\Gamma_1] \leq \mathbb{E}\int_{ t}^T e^{\lambda s + \mu A_s} |Y_s^{\frac{1}{n}} - Y_s^{\frac{1}{m}}|^2\mathrm{d}s+\frac{C}{N^2}\tilde{\mathcal{M}(T)}+\frac{C}{N^2}\mathcal{M}^2(T)+\frac{C}{N^2}\mathcal{M}(T).
    \end{equation*}
	Other estimates are as follows:
     \begin{equation*}
		\begin{aligned}
			&\Gamma_2+\Gamma_4\leq 2\int_{ t}^T e^{\lambda s + \mu A_s} |Y_s^{\frac{1}{n}} - Y_s^{\frac{1}{m}}|^2\mathrm{d}s\\
            &\quad+\int_{ t}^T e^{\lambda s + \mu A_s} |F_n(s, Y_s^{\frac{1}{n}}, Z_s^{\frac{1}{n}}, \mathcal{L}(Y_s^{\frac{1}{n}}, Z_s^{\frac{1}{n}})) - F(s, Y_s^{\frac{1}{n}}, Z_s^{\frac{1}{n}}, \mathcal{L}(Y_s^{\frac{1}{n}}, Z_s^{\frac{1}{n}}))|^21_{\bar{D}_{n,m}^N}\mathrm{d}s\\
            &\quad +\int_{ t}^T e^{\lambda s + \mu A_s} |F(s, Y_s^{\frac{1}{m}}, Z_s^{\frac{1}{m}}, \mathcal{L}(Y_s^{\frac{1}{m}}, Z_s^{\frac{1}{m}})) - F_m(s, Y_s^{\frac{1}{m}}, Z_s^{\frac{1}{m}}, \mathcal{L}(Y_s^{\frac{1}{m}}, Z_s^{\frac{1}{m}}))|^21_{\bar{D}_{n,m}^N}\mathrm{d}s ,\\
		\end{aligned}
	\end{equation*}
	\begin{equation*}
		\begin{aligned}
			&\Gamma_5\leq (K^2+1)\int_{ t}^T e^{\lambda s + \mu A_s}|Y_s^{\frac{1}{n}} - Y_s^{\frac{1}{m}}|^2\mathrm{d}A_s ,\\
		\end{aligned}
	\end{equation*}
    \begin{equation*}
		\begin{aligned}
			&\Gamma_3\leq \int_{ t}^T e^{\lambda s + \mu A_s} \big[(\frac{1}{8}+24r_N^2)|Y_s^{\frac{1}{n}} - Y_s^{\frac{1}{m}}|^2+\frac{1}{8}\|Z_s^{\frac{1}{n}}-Z_s^{\frac{1}{m}}\|^2\big] 1_{\bar{D}_{n,m}^N}\mathrm{d}s\\
			&\quad+\frac{1}{4}\int_{ t}^T e^{\lambda s + \mu A_s} \big[\mathbb{E}|Y_s^{\frac{1}{n}} - Y_s^{\frac{1}{m}}|^2+\mathbb{E}\|Z_s^{\frac{1}{n}}-Z_s^{\frac{1}{m}}\|^2\big] 1_{\bar{D}_{n,m}^N}\mathrm{d}s .\\
		\end{aligned}
    \end{equation*}
    Using the similar proof to Proposition 2, we have
    \begin{equation*}
		\begin{aligned}
			&e^{\lambda t +\mu A_{t }} |Y_{t }^{\frac{1}{n}} - Y_{t }^{\frac{1}{m}}|^2 
			+ \int_{ t}^T e^{\lambda s + \mu A_s} |Y_s^{\frac{1}{n}} - Y_s^{\frac{1}{m}}|^2 (\lambda \mathrm{d}s + \mu \mathrm{d}A_s) +\int_{ t}^T e^{\lambda s + \mu A_s} \|Z_s^{\frac{1}{n}} - Z_s^{\frac{1}{m}}\|^2 \mathrm{d}s  \\
			&\quad + 2 \int_{ t}^T e^{\lambda s + \mu A_s} \langle Y_s^{\frac{1}{n}} - Y_s^{\frac{1}{m}}, \nabla \varphi_n (Y_s^{\frac{1}{n}}) - \nabla \varphi_m (Y_s^{\frac{1}{m}}) \rangle \mathrm{d}s \\&\quad + 2 \int_{ t}^T e^{\lambda s + \mu A_s} \langle Y_s^{\frac{1}{n}} - Y_s^{\frac{1}{m}}, \nabla \psi_n (Y_s^{\frac{1}{n}}) - \nabla \psi_m (Y_s^{\frac{1}{m}}) \rangle \mathrm{d}A_s \\
            & \leq (\frac{25}{8}+24r_N^2)\int_{ t}^T e^{\lambda s + \mu A_s} |Y_s^{\frac{1}{n}} - Y_s^{\frac{1}{m}}|^2\mathrm{d}s +\frac{1}{8}\int_{ t}^T e^{\lambda s + \mu A_s} \|Z_s^{\frac{1}{n}}-Z_s^{\frac{1}{m}}\|^2 1_{\bar{D}_{n,m}^N}\mathrm{d}s\\
			&\quad+\frac{1}{4}\int_{ t}^T e^{\lambda s + \mu A_s} \big[\mathbb{E}|Y_s^{\frac{1}{n}} - Y_s^{\frac{1}{m}}|^2+\mathbb{E}\|Z_s^{\frac{1}{n}}-Z_s^{\frac{1}{m}}\|^2\big] 1_{\bar{D}_{n,m}^N}\mathrm{d}s\\
            &\quad+20L^2\big[\int_{ t}^T e^{\lambda s + \mu A_s} \eta_s^4\mathrm{d}s \big]^{\frac{1}{2}}\big[\int_{ t}^T e^{\lambda s + \mu A_s}(\dfrac{|Y_s^{\frac{1}{n}}|^2+|Y_s^{\frac{1}{m}}|^2+|Z_s^{\frac{1}{n}}|^2+|Z_s^{\frac{1}{m}}|^2}{N^2})^21_{D_{n,m}^N}\mathrm{d}s\big]^{\frac{1}{2}}\\
            &\quad +C\mathcal{M}(T)\mathbb{E}\int_{ t}^T e^{\lambda s + \mu A_s}\dfrac{|Y_s^{\frac{1}{n}}|^2+|Y_s^{\frac{1}{m}}|^2+|Z_s^{\frac{1}{n}}|^2+|Z_s^{\frac{1}{m}}|^2}{N^2}1_{D_{n,m}^N}\mathrm{d}s\\
            &\quad +\frac{40L^2}{N^2}\mathbb{E}\int_{ t}^T e^{\lambda s + \mu A_s}|Y_s^{\frac{1}{n}}|^{2}+|Y_s^{\frac{1}{m}}|^{2}+\|Z_s^{\frac{1}{n}}\|^{2}+\|Z_s^{\frac{1}{m}}\|^{2}1_{D_{n,m}^N}\mathrm{d}s\\
            &\quad +\frac{10L^2}{N^2}\mathbb{E}\int_{ t}^T e^{\lambda s + \mu A_s}\big(|Y_s^{\frac{1}{n}}|^{2}+|Y_s^{\frac{1}{m}}|^{2}+\|Z_s^{\frac{1}{n}}\|^{2}+\|Z_s^{\frac{1}{m}}\|^{2}\big)^21_{D_{n,m}^N}\mathrm{d}s \\
            &\quad+\int_{ t}^T e^{\lambda s + \mu A_s} |F_n(s, Y_s^{\frac{1}{n}}, Z_s^{\frac{1}{n}}, \mathcal{L}(Y_s^{\frac{1}{n}}, Z_s^{\frac{1}{n}})) - F(s, Y_s^{\frac{1}{n}}, Z_s^{\frac{1}{n}}, \mathcal{L}(Y_s^{\frac{1}{n}}, Z_s^{\frac{1}{n}}))|^21_{\bar{D}_{n,m}^N}\mathrm{d}s\\
            &\quad +\int_{ t}^T e^{\lambda s + \mu A_s} |F(s, Y_s^{\frac{1}{m}}, Z_s^{\frac{1}{m}}, \mathcal{L}(Y_s^{\frac{1}{m}}, Z_s^{\frac{1}{m}})) - F_m(s, Y_s^{\frac{1}{m}}, Z_s^{\frac{1}{m}}, \mathcal{L}(Y_s^{\frac{1}{m}}, Z_s^{\frac{1}{m}}))|^21_{\bar{D}_{n,m}^N}\mathrm{d}s\\
            &\quad +(K^2+1)\int_{ t}^T e^{\lambda s + \mu A_s}|Y_s^{\frac{1}{n}} - Y_s^{\frac{1}{m}}|^2\mathrm{d}A_s\\
             &\quad -2 \int_{ t}^T e^{\lambda s + \mu A_s} \langle \langle Y_s^{\frac{1}{n}} - Y_s^{\frac{1}{m}},(Z_s^{\frac{1}{n}} - Z_s^{\frac{1}{m}})dW_s\rangle.\\
        \end{aligned}
	\end{equation*}
    By taking the expectation on both sides of the previous estimate at $t=0$, and considering the constraints on $\lambda, \mu$ as well as Lemma 8, it follows that
    \begin{equation}
        \begin{aligned}
            &\mathbb{E} \int_{0}^T e^{\lambda s + \mu A_s} \big(|Y_s^{\frac{1}{n}} - Y_s^{\frac{1}{m}}|^2(\mathrm{d}s+\mathrm{d}A_s)+\|Z_s^{\frac{1}{n}} - Z_s^{\frac{1}{m}}\|^2 \mathrm{d}s \big)\\
            &\leq \frac{8}{5}e^{\frac{192Tr_N^2}{5}}\big\{\rho_N(F_n - F) +\rho_N(F_m - F)+\frac{C}{N^2}\tilde{\mathcal{M}(T)}+\frac{C}{N^2}\mathcal{M}^2(T)+\frac{C}{N^2}\mathcal{M}(T)\big\}.\quad
        \end{aligned}
	\end{equation}
    Taking the supremum over time and applying the BDG inequality, we obtain
    \begin{equation}
		\begin{aligned}
			&\mathbb{E}\sup_{0 \leq s \leq T}e^{\lambda s +\mu A_s} |Y_s^{\frac{1}{n}} - Y_s^{\frac{1}{m}}|^2\\
			&\leq [\frac{4}{3}+(C+\frac{256}{5}r_N^2)e^{\frac{192Tr_N^2}{5}}]\big\{\rho_N(F_n - F) +\rho_N(F_m - F)+\frac{C}{N^2}\tilde{\mathcal{M}(T)}+\frac{C}{N^2}\mathcal{M}^2(T)+\frac{C}{N^2}\mathcal{M}(T)\big\}.
		\end{aligned}
	\end{equation}
    Therefore,
    \begin{equation}
		\begin{aligned}
			&\mathbb{E}\sup_{0 \leq s \leq T}e^{\lambda s +\mu A_s} |Y_s^{\frac{1}{n}} - Y_s^{\frac{1}{m}}|^2+\mathbb{E} \int_{0}^T e^{\lambda s + \mu A_s} \big(|Y_s^{\frac{1}{n}} - Y_s^{\frac{1}{m}}|^2(\mathrm{d}s+\mathrm{d}A_s)+\|Z_s^{\frac{1}{n}} - Z_s^{\frac{1}{m}}\|^2 \mathrm{d}s \big) \\
			&\leq [\frac{4}{3}+(C+\frac{256}{5}r_N^2)e^{\frac{192Tr_N^2}{5}}]\big\{\rho_N(F_n - F) +\rho_N(F_m - F)+\frac{C}{N^2}\tilde{\mathcal{M}(T)}+\frac{C}{N^2}\mathcal{M}^2(T)+\frac{C}{N^2}\mathcal{M}(T)\big\}\\
            &\leq [\frac{4}{3}+(C+\frac{256}{5}r_N^2)e^{\frac{192Tr_N^2}{5}}]\big\{\rho_N(F_n - F) +\rho_N(F_m - F)+\frac{C}{N^{2(1-\alpha)}}[\tilde{\mathcal{M}(T)}+\mathcal{M}^2(T)+\mathcal{M}(T)]\big\}.
		\end{aligned}
	\end{equation}
    If $r_N^2\leq\dfrac{5(1-\alpha)}{96T}logN$, then by passing to the limit in $n,m,N$, we deduce that $(Y^{\frac{1}{n}}, Z^{\frac{1}{n}})$ forms a Cauchy sequence in the Banach space $(\mathcal{S}_k^{\lambda, \mu} \cap \mathcal{H}_k^{\lambda, \mu} \cap \widetilde{\mathcal{H}}_k^{\lambda, \mu}) \times \mathcal{H}_{k \times d}^{\lambda, \mu}$. Hence, there exists $(Y, Z)$ such that
	\begin{equation}
		\begin{aligned}
			&\lim_{n \to \infty}\mathbb{E} \int_0^T e^{\lambda s + \mu A_s}  |Y_s^{\frac{1}{n}} - Y_s|^2\mathrm{d}s =0 ,\quad \lim_{n \to \infty}\mathbb{E} \int_0^T e^{\lambda s + \mu A_s}  |Y_s^{\frac{1}{n}} - Y_s|^2\mathrm{d}A_s =0,\\
			&\lim_{n \to \infty}\mathbb{E} \int_0^T e^{\lambda s + \mu A_s}  \|Z_s^{\frac{1}{n}} - Z_s\|^2\mathrm{d}s =0  ,\quad\lim_{n \to \infty} \mathbb{E} \sup_{0 \leq s \leq T } e^{\lambda s + \mu A_s} |Y_s^{\frac{1}{n}} - Y_s|^2=0 .
		\end{aligned}
	\end{equation}
    
    For any given $N$, we define $H_n^N:=\{(s,w): |Y_s^{\frac{1}{n}}|^2+|Y_s|^2+|Z_s^{\frac{1}{n}}|^2+|Z_s|^2\geq N^2\}$ and $\bar{H}_n^N=\Omega-H_n^N$.By an argument analogous to the estimations of $\Gamma_1$ and $\Gamma_3$, we have
	\begin{equation*}
		\begin{aligned}
			&\mathbb{E}\int_0^T e^{\lambda s + \mu A_s}  |F_n(s, Y_s^{\frac{1}{n}}, Z_s^{\frac{1}{n}},\mathcal{L}(Y_s^{\frac{1}{n}}, Z_s^{\frac{1}{n}}))-F(s, Y_s, Z_s,\mathcal{L}(Y_s, Z_s))|^2 \mathrm{d}s\\
			&\leq \frac{C}{N^2}\tilde{\mathcal{M}(T)}+\frac{C}{N^2}\mathcal{M}^2(T)+\frac{C}{N^2}\mathcal{M}(T)+2\rho_N(F_n-F)\\
            &\quad +16L^2\mathbb{E}\int_0^T e^{\lambda s + \mu A_s}\|Z_s^{\frac{1}{n}}-Z_s\|^2\mathrm{d}s +16L^2\mathbb{E}\int_0^T e^{\lambda s + \mu A_s}|Y_s^{\frac{1}{n}}-Y_s|^2\mathrm{d}s .\quad\\
		\end{aligned}
	\end{equation*}
    Passing to the limit in $n$ and $N$ , we obtain
	\begin{equation}
		\mathbb{E}\int_0^T e^{\lambda s + \mu A_s}  |F_n(s, Y_s^{\frac{1}{n}}, Z_s^{\frac{1}{n}},\mathcal{L}(Y_s^{\frac{1}{n}}, Z_s^{\frac{1}{n}}))-F(s, Y_s, Z_s,\mathcal{L}(Y_s, Z_s))|^2 \mathrm{d}s \to 0 .
	\end{equation}

    Applying the same approach as in section 2.2.2 to the results above, we establish the existence.
    
    The proof of uniqueness follows the same method as in Proposition 2 and is therefore omitted.
	
	\newpage
	\bibliographystyle{plain}
	\bibliography{referencelist}
	
\end{document}